\documentclass[12pt]{article}

\usepackage[dvipsnames]{xcolor}
\usepackage{enumerate,enumitem}
\usepackage{amsmath,amssymb,amsfonts,mathrsfs, amsthm,mathtools, bm, bbm, dsfont, mathrsfs, amsthm}
\usepackage{graphicx, epstopdf}

\usepackage{fullpage}

\usepackage{amsthm}
\allowdisplaybreaks

\usepackage{amsmath,amssymb,amsfonts}

\usepackage{graphicx}
\usepackage{multicol}
\usepackage{tikz}
\usepackage{amsfonts}
\usepackage{scrextend}
\usepackage{amsfonts}
\usepackage{amsmath}
\usepackage{amssymb}
\usepackage{enumitem}
\usepackage{mathtools}
\usepackage{breqn}
\usepackage{dsfont}
\usepackage{xcolor}
\usepackage[colorlinks]{hyperref}

\DeclareMathOperator*{\argminA}{arg\,min} 

\usetikzlibrary{fit}
\hypersetup{
colorlinks=true,       
linkcolor=blue,          
citecolor=blue,        
filecolor=blue,      
urlcolor=blue 
}

\DeclareMathOperator{\var}{Var}

\newcommand{\X}{\mathbb{X}}
\newcommand{\U}{\mathbb{U}}
\newcommand{\Pc}{\mathcal{P}}
\newcommand{\LEXN}{L_{\sf EX}^{N}}
\newcommand{\LN}{L^{N}}
\newcommand{\PN}{$\mathcal{P}^{N}$}
\newcommand{\PNT}{$\mathcal{P}^{N}_T$}

\newcommand{\PIN}{$\mathcal{P}^{\infty}$}
\newcommand{\PIT}{$\mathcal{P}^{\infty}_T$}

\newcommand{\LCONS}{L_{\sf CO,SYM}^{N}}
\newcommand{\LPRNS}{L_{\sf PR,SYM}^{N}}

\newcommand{\LCOS}{L_{\sf CO,SYM}}
\newcommand{\LPRS}{L_{\sf PR,SYM}}
\newcommand{\LEX}{L_{\sf EX}}

\newcommand{\PPN}{\Pi^{N}}
\newcommand{\PM}{\Pi^{N}_{\sf M}}
\newcommand{\PS}{\Pi^{N}_{\sf S}}

\newcommand{\PCONS}{\Pi_{\sf CO,SYM}^{N}}
\newcommand{\PPRNS}{\Pi_{\sf PR,SYM}^{N}}

\newcommand{\PCOS}{\Pi_{\sf CO,SYM}}
\newcommand{\PPRS}{\Pi_{\sf PR,SYM}}
\newcommand{\PEX}{\Pi_{\sf EX}}
\newcommand{\PEXN}{\Pi_{\sf EX}^{N}}
\newcommand{\PEXNM}{\Pi_{\sf EX,M}^{N}}
\newcommand{\PEXNS}{\Pi_{\sf EX,S}^{N}}
\newcommand{\PEXM}{\Pi_{\sf EX,M}}
\newcommand{\PEXS}{\Pi_{\sf EX,S}}

\newcommand{\EMX}{\mathcal{P}^N_{\sf{E}}(\mathbb{X})}
\newcommand{\EM}{\mathcal{P}^N_{\sf{E}}(\mathbb{X}\times \mathbb{U})}

\newtheorem{theorem}{Theorem}

\theoremstyle{definition}
\newtheorem{definition}{Definition}
\newtheorem{lemma}{Lemma}
\newtheorem{assumption}{Assumption}
\newtheorem{example}{Example}
\newtheorem{corollary}{Corollary}
\theoremstyle{remark}

\date{}
\begin{document}

\title{Optimality of Symmetric Independent  Policies under Decentralized Mean-Field Information Sharing for Stochastic Teams and Equivalence with McKean-Vlasov Control of a Representative Agent}

\author{Sina Sanjari\thanks{Sina Sanjari is with the Department of Mathematics and Computer Science, Royal Military College of Canada, Canada; {\tt\small sanjari@rmc.ca}} \quad Naci Saldi\thanks{Naci Saldi is with the Department of Mathematics, Bilkent University, Ankara, Turkey; {\tt\small naci.saldi@bilkent.edu.tr}} \quad Serdar Y\"uksel\thanks{Serdar Y\"uksel is with the Department of Mathematics and Statistics, Queen's University, Canada; {\tt\small yuksel@queensu.ca}}}
\maketitle

\begin{abstract}
We study a class of stochastic exchangeable teams with a finite number of decision makers (DMs) as well as their mean-field limits with infinitely many DMs. In the finite population regime, we study exchangeable teams under the centralized information structure. For the infinite population setting, we study both the centralized information structure and the decentralized mean-field information sharing structure. The paper makes the following main contributions: i) For finite population exchangeable teams, we establish the existence of an optimal policy that is exchangeable (permutation invariant) and Markovian; ii) As our main result in the paper, we show that a sequence of exchangeable optimal policies for finite population settings (which satisfies a measure valued MDP formulation due to B{\"a}uerle) converges to a decentralized symmetric (identical) and conditionally independent (given the mean-field) policy for the infinite population problem, which is then globally optimal under both the centralized information structure as well as the mean-field sharing information structure. (iii) This result establishes existence of a symmetric, independent, decentralized optimal randomized policy for the infinite population problem and proves the optimality of the limiting measure-valued MDP for the representative DM. Our paper thus establishes the relation between controlled McKean-Vlasov dynamics and the optimal infinite population decentralized stochastic control problem (without an apriori restriction of symmetry in policies of individual agents), for the first time, to our knowledge (beyond several special cases). We also establish near optimality of a numerical method for solving this problem. iv) Finally, we show that symmetric, independent, decentralized optimal randomized policies are approximately optimal for the corresponding finite-population team with a large number of DMs under the centralized information structure.
\end{abstract}

\section{Introduction.}
Decentralized stochastic control, or teams, involve a collection of decision makers (DMs) acting together to optimize a common cost but not necessarily sharing all their available information. Each DM has partial access to global information determined by the information structure, as in \cite{wit75,YukselBasarBook24}. This differs from stochastic games, where DMs may possess conflicting cost criteria, probabilistic models, or priors on the system. Stochastic teams (see e.g., \cite{ho1980team, CDCTutorial}) extend classical single DM stochastic control problems, finding applications in various fields such as networked control as in \cite{ho1980team, hespanha2007survey, YukselBasarBook}, communication networks as in \cite{hespanha2007survey}, cooperative systems as in \cite{mar55, Radner,carrillo2022controlling}, and large sensor networks as in \cite{tsitsiklis1988decentralized}.

In the context of stochastic games with many DMs, a significant body of work revolves around mean-field (MF) games: these games serve as limiting models of weakly-interacting symmetric finite DM games (see e.g., \cite{CainesMeanField2,LyonsMeanField,carmona2020applications}). Several studies, such as \cite{LyonsMeanField,bardi2019non, carmona2016mean,lacker2015mean,saldi2018markov}, have established the existence of a symmetric Nash equilibrium for MF games. On one front, the MF approach devises policies such that Nash equilibria for MF games are demonstrated to be approximately Nash equilibria for the corresponding (pre-limit) games with large numbers of DMs; see e.g., \cite{CainesMeanField3, saldi2018markov, carmona2018probabilistic, cecchin2017probabilistic}. Conversely, the limits of symmetric Nash equilibrium policies of finite DM games are shown to converge to a Nash equilibrium of the corresponding MF game as the number of DMs tends to infinity; see e.g., \cite{fischer2017connection, lacker2018convergence, LyonsMeanField}. 

Here, we study discrete-time stochastic MF teams under global information (or centralized information) as well as local MF information sharing, where agents have access to the empirical measures of states in addition to their local states and actions. Our primary objective is to establish the existence and structural characteristics (such as symmetry or exchangeability) of a globally optimal solution for teams with a large but finite number of decision makers (DMs) under the global information structure, as well as their infinite population limit MF teams with decentralized MF information sharing. These existence and structural results are significant in the development of computational, approximation, and learning methodologies.  Specifically, we leverage these results to demonstrate the optimality of a single representative DM's Markov decision problem (MDP) for the infinite population team and its approximation for finite but large population teams. Our approach differs from the methods employed in MF games. This distinction primarily arises because Nash equilibria or person-by-person optimal solutions, within the context of teams, often correspond to locally optimal solutions. Consequently, the existence and structural outcomes for games may yield inconclusive insights regarding global optimality for teams without uniqueness (e.g., see \cite{LyonsMeanField, hajek2019non, bayraktar2020non,bardi2019non,delarue2020selection} for non-uniqueness of Nash equilibrium). 

{In \cite{bauerle2023mean}, an equivalent MDP formulation is characterized for finite population MF teams with the centralized information structure. Furthermore, in  \cite[Theorem 4.6]{bauerle2023mean}, it has been shown that the limsup of the sequence value functions for the finite population MDP is equal to the value function of a limiting MDP as the number of DMs goes to infinity and this certifies that the value function of the problem is attained as the solution of the measure-valued MDP. \cite{bauerle2023mean} did not establish that the optimal value function can be achieved by independently randomized symmetric policies, except for the case when dynamics are decoupled (that is, when the flow is i.i.d); see \cite[Remark 4.7b]{bauerle2023mean}.} Finite approximation models for finite population MDPs have also been studied in \cite{Ali-2023} using the MDP connections in \cite{bauerle2023mean}, by quantizing the state space. In addition, considering symmetric policies for the infinite population MF problem, a finite approximation MDP is obtained for the limiting MDP \cite{Ali-2023} by discretizing the space of probability measures on the states of a representative agent. Related to \cite{bauerle2023mean}, in \cite{carmona2023model}, dynamic programming equations have been established for the limiting MF control (McKean-Vlasov) problem with a representative agent using a lifted measure-valued MDP. Unlike \cite{bauerle2023mean}, \cite{carmona2023model} did not study the relation with the finite or infinite population problem studied in \cite{bauerle2023mean}, and the connection between the finite population problem and the MF control problem has not been investigated in \cite{carmona2023model}. In other words, it is assumed that apriori policies are symmetric in the infinite population limit (see \cite[Section 2]{carmona2023model}). Additionally, in \cite{motte2022mean}, the connection between the finite population problem and the limiting MDP for the McKean-Vlasov MDP has been established under the assumption that the policies are open-loop, symmetric, and decentralized. Furthermore, in \cite{Mahajan-CDC14}, the MF sharing information structure has been considered for finite population MF teams under symmetric policies, where dynamic programming equations have been obtained using the common information approach. {In this paper, building on the MDP formulation in \cite{bauerle2023mean}, we utilize a proof technique in our prior work in \cite[Lemma 1]{SSYdefinetti2020} (constructing an exchangeable policy that attains the same performance as any given policy) to demonstrate that the solution of the measure-valued MDP in \cite{bauerle2023mean} can be realized as an exchangeable policy (possibly asymmetric) for the original team problem under the centralized information structure with finite population DMs.} Subsequently, for the infinite population team, we further leverage exchangeability and a de Finetti-type representation theorem from \cite{SSYdefinetti2020} to establish the existence of a symmetric and decentralized optimal solution with MF information sharing. Consequently, we establish that the solution of the representative-DM MDP can be realized as a symmetric decentralized policy that is optimal for the original team problem with infinite population DMs. To our knowledge, this is the first result on the optimality of symmetric policies leading to the optimality of a representative DM-based optimization, thus connecting existing results on measure-valued MDPs and McKean-Vlasov dynamics with decentralized stochastic control.

In prior works \cite{sanjari2018optimal, sanjari2019optimal, SSYdefinetti2020}, we studied discrete-time MF teams under strict decentralized information structure (without MF sharing), which do not allow for a Markovian formulation and primarily build on change of measure arguments. In \cite{sanjari2019optimal}, we studied convex MF teams with decoupled dynamics under a decentralized information structure. Leveraging convexity and a lower semi-continuity argument, we established the existence of a symmetric optimal solution for both finite and infinite populations. In \cite{SSYdefinetti2020}, we extended the analysis to non-convex MF teams under the same decentralized setting, employing exchangeability and a de Finetti theorem to demonstrate the existence of an optimal solution that is asymptotically symmetric and independent. In contrast to \cite{sanjari2019optimal,SSYdefinetti2020}, which focused on decentralized information structures, the present work begins with the centralized information structure. It then studies the limiting MF problem under the decentralized MF-sharing information structure. While we also rely on exchangeability and a de Finetti-type argument as in \cite{SSYdefinetti2020}, the argument is fundamentally different: (i) the centralized nature of the information structure prevents us from applying the change-of-measure technique and the compactification and convergence analysis for the policy space used in \cite{SSYdefinetti2020}, and (ii) critically we show that the limiting policies are conditionally independent given the empirical distribution, and this conditional independence is attainable under the MF sharing pattern. To this end, we adopt the measure-valued MDP formulation developed in \cite{bauerle2023mean}, along with its value iteration and strategic measures defined over empirical state distributions. Furthermore, starting with a centralized information structure enables an explicit connection to the MF sharing problem and facilitates the analysis of value iteration in lifted MDPs—directions that were not explored in \cite{SSYdefinetti2020}. Additionally, in \cite{sanjari2024nash}, which builds on \cite{SSYdefinetti2020}, team vs team has been studied in a game formulation under a strict decentralized information structure. In particular, the properties of a Nash equilibrium that are globally optimal within DMs of each team have been studied. In this paper, we study a single team formulation under centralized and MF sharing information structures.

A notable set of results related to MF teams is that on social optimal control problems, MF optimal control problems, and their limiting problems of McKean-Vlasov optimal control (see e.g., \cite{motte2022mean,carmona2023model,cecchin2021finite,carmonaForwardbackward,pham2017dynamic,bensoussan2015master,bayraktar2018randomized,djete2022mckean,lacker2017limit,motte2023quantitative,fornasier2019mean,albi2017mean,albi2022moment,lauriere2016dynamic,achdou2020mean}), with particular emphasis on the linear quadratic Gaussian model (see e.g., \cite{huang2012social,carmona2018probabilistic, caines2018peter, pham2016discrete,elliott2013discrete,ni2015discrete,toumi2024mean,arabneydi2015team}).  Within this framework, characterizing the bounds for the discrepancy between centralized social performance and decentralized controllers has been a central focus under open-loop, relaxed, and Markovian feedback policies. Results on McKean-Vlasov optimal control problems, or more generally MF-type control problems, assume symmetry in policies among DMs to make a connection between the infinite population problem and the optimal control problem of the representative player. 

Our results are also closely related to \cite{jackson2023approximately,cardaliaguet2023algebraic} in the context of MF settings, where non-asymptotic bounds have been derived for the gap between the optimal performances under centralized and strictly decentralized information structures. {In particular, in \cite{jackson2023approximately}, a continuous-time MF problem has been studied under both fully centralized and decentralized information structures. Under a convexity condition,  \cite{jackson2023approximately} obtained quantitative bounds between centralized and decentralized value functions for a particular class of problems with decoupled dynamics. Focusing on a specific class of teams allows the derivation of HJB equations for both the centralized problem and decentralized problem by lifting the states to their laws (see \cite[equations (3.9) and (3.13)]{jackson2023approximately}).  For the open-loop decentralized information structure, drawing a connection to the HJB equation for person-by-person optimality allows the maximum principle to be applied to obtain the necessary conditions for optimality. Using the convexity of the problem and uniqueness in person-by-person optimality, a quantitative bound is established using Poincar\`{e} and transport inequalities (see \cite[Theorem 4.5]{jackson2023approximately}). We note that symmetry of an optimal decentralized policy is utilized to establish an upper bound between the decentralized optimal cost and that of the MF limit. The bound follows from the convergence of empirical measures of i.i.d. random variables. We note that the symmetry of an optimal policy has not been explicitly established, but due to convexity and using uniqueness and person-by-person optimality to characterize the optimal solution, it can be concluded for this specific team problem. This is, however, not true for a more general setting of non-convex exchangeable MF teams with finite $N$-agents studied in this paper. Here we consider a general (possibly non-convex) setting of discrete-time MF teams with coupled dynamics. In contrast to some of the above results, we do not assume symmetry a priori.}

The paper makes the following main contributions. 

\begin{enumerate}
\item In Theorem \ref{the:2}, we show that a solution of an equivalent measure-valued MDP can be realized as an exchangeable Markov policy (also stationary for the infinite horizon problem), but not necessarily symmetric (i.e., identical across DMs) under the centralized information structure. 
\item In Theorem \ref{the:3}, we demonstrate the convergence of a sequence of exchangeable optimal policies, along with the corresponding value iteration, from a finite population setting to a decentralized randomized policy for the infinite-population problem that is conditionally symmetric and conditionally independent, together with its corresponding value iteration. This convergence establishes the existence of a globally optimal solution for the infinite population MF team and shows that an optimal solution is independently randomized, symmetric, and decentralized with MF information sharing. Building on this symmetry and independence (over agents conditioned on MF information) of an optimal policy, we establish the optimality of the MDP for the representative DM in Corollary \ref{Cor:MV} and Theorem \ref{the:veri}. Our paper thus establishes the relation between the controlled McKean-Vlasov dynamics and the optimal infinite population decentralized stochastic control problem (without an apriori restriction of symmetry in policies of individual agents), for the first time, to our knowledge.

\item In Theorem \ref{the:4}, we show that a symmetric optimal solution with the decentralized MF-sharing information structure for the MF team is near optimal for the finite-population team with a large number of DMs. In Theorems  \ref{the:veri} and \ref{the:Approx-MK}, we provide a Bellman recursion and, via its weak Feller and compact MDP formulation, establish near optimality of a finite model. 
\end{enumerate}

A summary of our results is depicted in Fig. \ref{fig:1}, and our proof technique is explained in Section \ref{sec:sum}.

\section{Exchangeable Teams and Their MF Limit.}

In this section, we first introduce a class of exchangeable teams with a finite number of decision makers (DMs), $N \in \mathbb{N}$, and then we introduce their infinite population MF teams.

\subsection{Finite population exchangeable teams.}

The dynamics of each DM$^{i}$
are given by
\begin{align}
x_{t+1}^{i}=f(x_t^i,u_t^i,\mu_t^N,w_t^i) \quad \forall\: t\geq 0\label{eq:dyn}
\end{align}
for $i\in \mathcal{N}:=\{1,\ldots,N\}$.
In the above dynamics, the process $\{w_{t}^{i}\}_{t\geq0}$ is an i.i.d. process  which is independent of i.i.d. random variables $x_{0}^{i}\sim P_0$ for all $i\in \mathcal{N}$. The state, action, and disturbance spaces for each DM are assumed to be subsets of appropriate dimensional Euclidean spaces and are denoted by $\X$, $\U$, and $\mathbb{W}$, respectively. The measure-valued process $\{\mu_t^N\}_{t\geq 0}$ is the empirical measure process on states, i.e., for every $t\geq0$ 
\begin{align*}
    \mu_t^N(\cdot):=\frac{1}{N}\sum_{i=1}^{N}\delta_{x_{t}^{i}}(\cdot),
\end{align*}
where $\delta_x(\cdot)$ is the Dirac delta measure at point $x$. Let $\mathcal{P}(\X)$ be the set of all probability measures on $\X$ endowed with the weak convergence topology.
Let $\EMX\subset \mathcal{P}(\X)$ be the set of all empirical measures on $\mathbb{X}$ with $N$ DMs, i.e.,
\begin{align*}
    \EMX := \left\{\frac{1}{N}\sum_{i=1}^{N}\delta_{x^{i}}(\cdot) \bigg| x^i\in \mathbb{X}, \quad \forall i\in \mathcal{N}\right\}.
\end{align*}

We note that $\mu_t^N\in \EMX$ for all $t\geq 0$. Throughout the paper, we use the notation $x^{1:N}_{t}$ to denote $(x^{1}_{t}, \ldots, x^{N}_{t})$ for all $t\geq 0$. 

Denote the set of all stochastic kernels on $\mathbb{A}$ given random variables taking values in $\mathbb{B}$ by $\mathcal{P}(\mathbb{A}|\mathbb{B})$. 
Following \eqref{eq:dyn}, random variables $x^{1:N}_{t+1}$ are jointly distributed according to 
\begin{align}\label{st-ker-1}
    x^{1:N}_{t+1} \sim \prod_{i=1}^{N}{\mathcal{T}}(\cdot|x_{t}^i,u_{t}^i, \mu_{t}^N),
\end{align}
where ${\mathcal{T}}$ is a transition kernel in $\mathcal{P}(\X|\X \times \U \times \mathcal{P}(\X))$ representing dynamics in \eqref{eq:dyn}. We emphasize that due to symmetry in dynamics among DMs, ${\mathcal{T}}$ is the same for all DMs, and also since dynamics are time-invariant, so is ${\mathcal{T}}$.

We consider the classical (or centralized) information structure for each DM, which is given by 
\begin{align}\label{eq:IS}
    I_{t}^{i}:=I_{t}=\{x_{0:t}^{1:N},u_{0:t-1}^{1:N}\}\quad \forall t\geq 0.
\end{align}
Under the centralized information structure, DMs have full access to the history of states and actions and the current state of all DMs. {This information structure, however, is restrictive, specifically for teams with large numbers of DMs, since it requires that all available information be shared among all DMs.} We discuss later the decentralized information structure with MF sharing (see \eqref{eq:IS-DEC}) and discuss the connection between its corresponding team solution and that of centralized information structure.

Our goal in this paper is to study both finite and infinite horizon problems with finite and infinite populations of DMs. We now define admissible deterministic policies for each DM under the centralized information structure. The set of deterministic policies for each DM$^i$ at time $t\geq 0$ is the set of measurable functions $\gamma_{t}^{i}$ that map the random variables in the information set $I_{t}$ to the DM$^{i}$'s action $u^{i}_{t}$. We denote the set of all the admissible policies at time $t\geq 0$ by $\Gamma^i_t$. The set of admissible policies $\pmb{\gamma^i_T}=(\gamma_{0}^i, \gamma_{1}^i, \ldots, \gamma_{T-1}^i)$ for each DM$^i$ over a finite horizon of length $T$ is denoted by $\pmb{\Gamma^i_T}:=\prod_{t=0}^{T-1}\Gamma^i_t$ and over infinite horizon $\pmb{\gamma^i}=(\gamma_{0}^i, \gamma_{1}^i, \ldots)$ is denoted by $\pmb{\Gamma^i}:=\prod_{t=0}^{\infty}\Gamma^i_t$. Throughout the paper, we use the bold notation for any variable $\pmb{a}$ to denote the infinite sequence $\pmb{a}:=(a_1,a_2 \ldots)$ over infinite horizon. Also, we use the one with subscript $T$ for the length $T$ sequence $\pmb{a_T}:=(a_1,\ldots,a_T)$.

We next define both finite-horizon and infinite-horizon discounted expected costs. The finite horizon discounted expected cost for any admissible policy $\pmb{\gamma^{1:N}_T}:=(\pmb{\gamma^{1}_T},\ldots, \pmb{\gamma^{N}_T})\in \prod_{i=1}^{N}\pmb{\Gamma^i_T}$ is given by
\begin{align}\label{expected-cost-N-T}
    J^N_T(\pmb{\gamma^{1:N}_T})=\mathbb{E}^{\pmb{\gamma^{1:N}_T}}\left[\frac{1}{N}\sum_{i=1}^{N}\sum_{t=0}^{T-1}\beta^tc(x_{t}^{i},u_{t}^{i},\mu_{t}^N)\right],
\end{align}
where the cost $c:\X \times \U \times \mathcal{P}(\X)\to \mathbb{R}_{+}$ is a Borel measurable function and the discount factor $\beta \in (0,1)$.  We call this team problem \emph{exchangeable} because the cost $$\widehat{c}_T\left(\pmb{x^{1:N}_{T}},\pmb{u^{1:N}_{T}}\right):=\frac{1}{N}\sum_{i=1}^{N}\sum_{t=0}^{T-1}\beta^tc(x_{t}^{i},u_{t}^{i},\mu_{t}^N)$$ satisfies the following exchangeability condition 
\begin{align*}
\widehat{c}_T\left(\pmb{x^{1:N}_{T}},\pmb{u^{1:N}_{T}}\right)=\widehat{c}_T\left(\pmb{x^{\sigma(1):\sigma(N)}_{T}},\pmb{u^{\sigma(1):\sigma(N)}_{T}}\right)
\end{align*}
for all permutations $\sigma$ of $\{1,\ldots, N\}$, where ${a^{\sigma(1):\sigma(N)}}:={a^{\sigma(1)}},\ldots, {a^{\sigma(N)}}$ for $a=\pmb{x_T}$ or $\pmb{u_T}$.

The infinite horizon discounted expected cost for any admissible policy $\pmb{\gamma^{1:N}}:=(\pmb{\gamma^{1}},\ldots, \pmb{\gamma^{N}})\in \prod_{i=1}^{N}\pmb{\Gamma^i}$ is given by
\begin{align}\label{expected-cost-N-inf}
    J^N(\pmb{\gamma^{1:N}})=\mathbb{E}^{\pmb{\gamma^{1:N}}}\left[\frac{1}{N}\sum_{i=1}^{N}\sum_{t=0}^{\infty}\beta^tc(x_{t}^{i},u_{t}^{i},\mu_{t}^N)\right].
\end{align}
We denote the finite horizon team problem by \PNT\ while the infinite horizon discounted team problem is denoted by \PN.
Next, we define the notion of global optimality for \PN\ (a similar definition holds for \PNT).

\begin{definition}\label{eq:gof}
For the stochastic team \PN, a policy $\pmb{{\gamma}^{1\star:N\star}}:=(\pmb{{\gamma}^{1\star}},\ldots, \pmb{{\gamma}^{N\star}})\in \prod_{i=1}^{N}\pmb{\Gamma^{i}}$\@ is
globally optimal for \PN\ if
\begin{equation*}
J^{N}(\pmb{{\gamma}^{1\star:N\star}})=\inf_{\pmb{{\gamma}^{1:N}}\in \prod_{i=1}^{N}\pmb{\Gamma^{i}}}
J^{N}(\pmb{{\gamma}^{1:N}}).
\end{equation*}  
\end{definition}

In the next subsection, we introduce the infinite population MF limit of \PNT\ and \PN\ by taking the limit of the number of DMs as $N\to \infty$.

\subsection{Infinite population MF teams.}

In this section, we define the limiting infinite population teams. For the infinite population finite horizon discounted MF team, the expected cost for any admissible policy $\underline{\pmb{\gamma_T}}:=(\pmb{\gamma^1_T},\pmb{\gamma^2_T}, \ldots)\in \prod_{i=1}^{\infty}\pmb{\Gamma^i_T}$ is given by
\begin{align}\label{expected-cost-inf-T}
    J(\underline{\pmb{\gamma_T}})=\limsup_{N\to \infty}\mathbb{E}^{\underline{\pmb{\gamma_T}}}\left[\frac{1}{N}\sum_{i=1}^{N}\sum_{t=0}^{T-1}\beta^t c(x_{t}^{i},u_{t}^{i},\mu_{t}^N)\right],
\end{align}
where the cost $c:\X \times \U \times \mathcal{P}(\X)\to \mathbb{R}_{+}$ is a Borel measurable function. Throughout the paper, we use $\underline{a}$ for any variable $a$ to denote the infinite sequence of variables $(a^{1}, a^2, \ldots)$ among DMs.

For the infinite population infinite horizon discounted MF team, the expected cost for any admissible policy $\underline{\pmb{\gamma}}:=(\pmb{\gamma^1},\pmb{\gamma^2}, \ldots)\in \prod_{i=1}^{\infty}\pmb{\Gamma^i}$ is given by
\begin{align}\label{expected-cost-inf-inf}
    J(\underline{\pmb{\gamma}})=\limsup_{N\to \infty}\mathbb{E}^{\underline{\pmb{\gamma}}}\left[\frac{1}{N}\sum_{i=1}^{N}\sum_{t=0}^{\infty}\beta^tc(x_{t}^{i},u_{t}^{i},\mu_{t}^N)\right],
\end{align}
where the cost $c:\X \times \U \times \mathcal{P}(\X)\to \mathbb{R}_{+}$ is a Borel measurable function. 

We denote the infinite population finite horizon MF team problem by \PIT\ while the infinite population infinite horizon discounted  MF team problem is denoted by \PIN.

Next, we define the notion of global optimality for \PIN\ (a similar definition holds for \PIT).

\begin{definition}\label{eq:gif}
For stochastic team \PIN, a policy $\underline{\pmb{{\gamma}^{\star}}}:=(\pmb{{\gamma}^{1\star}},\pmb{{\gamma}^{2\star}},\ldots)\in \prod_{i=1}^{\infty}\pmb{\Gamma^{i}}$\@ is
globally optimal for \PIN\ if
\begin{equation}
J(\underline{\pmb{{\gamma}^{\star}}})=\inf_{\underline{\pmb{{\gamma}}}\in \prod_{i=1}^{\infty}\pmb{\Gamma^{i}}}
J(\underline{\pmb{{\gamma}}}).
\end{equation}  
\end{definition}

Our main results use the following set of assumptions.

\begin{assumption}\label{assump:1}
The following three statements hold:
\begin{itemize}
    \item [(i)] $\U$ and $\X$ are compact;
    \item [(ii)] $c:\X \times \U \times \mathcal{P}(\mathbb{X}) \to \mathbb{R}_{+}$ is jointly continuous  in all arguments;
    \item [(iii)] $f(\cdot,\cdot,\cdot, w):\X \times \U \times \mathcal{P}(\X) \to \X$ in \eqref{eq:dyn} is jointly continuous for all $w\in \mathbb{W}$.
\end{itemize}
\end{assumption}

We note that compactness of $\X$ can be relaxed as our analysis only requires the sequence of probability measures on $\X$ to be weakly compact. Our main results require randomization in policies for DMs within the team problems. We first provide a high-level description of our main results.

\section{Summary and discussion on main results.}\label{sec:sum}

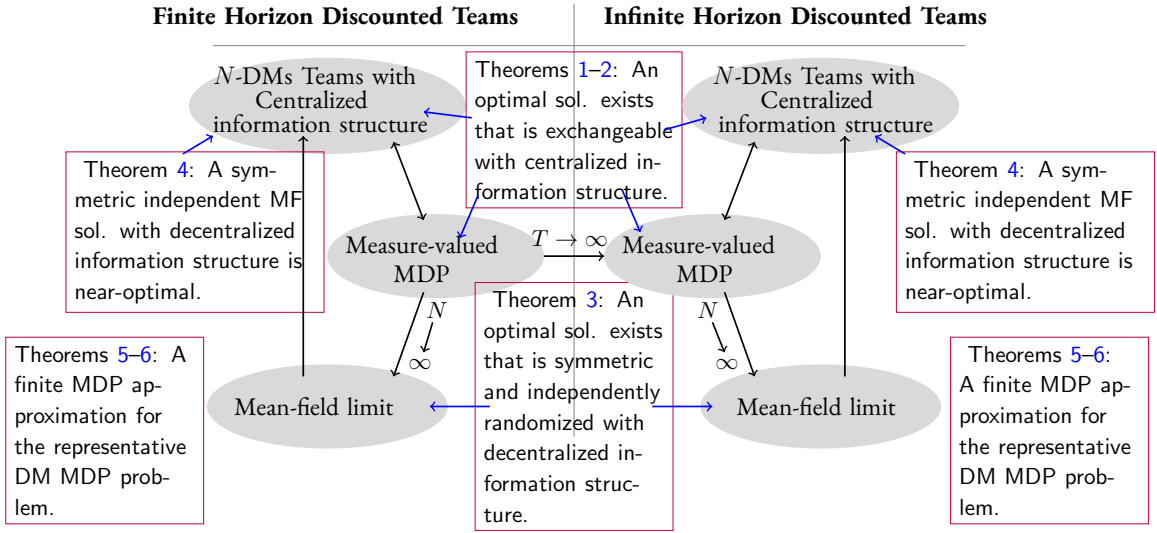
\begin{figure*}[h]
\centering
 \scalebox{0.8}{
    \begin{tikzpicture}[very thin]
    \fill[black!15] (-0.6,0) ellipse (2.3 and 0.8);
    \fill[black!15] (-0.6,-5) ellipse (2 and 0.7);
    \fill[black!15] (1.2,-2.5) ellipse (1.8 and 0.7);

    \draw[gray] (3.5,1.7) -- (3.5,-5.5);
    \draw[gray] (-2.5,1) -- (10,1);
    \node[text width=8cm] at (0.5,1.5) 
    {\bf Finite Horizon Discounted Teams};
    \node[text width=4cm] at (-0.5,0.5) 
    {$N$-DMs Teams with};
    \node[text width=4cm] at (0.2,0.1) 
    {Centralized};
    \node[text width=4cm] at (-0.5,-0.3) 
    {information structure};
    \node[text width=4cm] at (-0.15,-5) 
    {Mean-field limit};
    \node[text width=3cm] at (1.2,-2.3) 
    {Measure-valued};
    \node[text width=3cm] at (2,-2.8) 
    {MDP};
    \node[black,text width=3cm] at (-4.3,-2.1) (the4)
    {\small{ \sf Theorem \ref{the:4}: A symmetric independent MF sol. with decentralized information structure is near-optimal.}} ;
    \node[black,text width=3cm] at (3.5,-5.3) (the3)
    {\small{ \sf Theorem \ref{the:3}: An optimal sol. exists that is symmetric and independently randomized with a decentralized information structure.}};
    \node [rectangle,purple,draw,fit=(the4),inner sep=0] {};
    \node [rectangle,purple,draw,fit=(the3),inner sep=0] {};

    \node[black,text width=0.1cm] at (1.1,-3.4) 
    {\small{$N$}};
    \node[black,text width=0.1cm] at (0.8,-4.3) 
    {\small{$\infty$}};
    \draw[black, thick, ->]  (1.15,-3.6) -- (1,-4.1);

    \draw[black, thick, ->]  (-1,-4.5) -- (-1,-0.5);
    \draw[black, thick, <->]  (0.5,-0.6) -- (1,-1.9);
    \draw[black, thick, <-]  (0.5,-4.5) -- (1,-3.1);
\node[black,text width=3cm] at (11.4,-2.1) (the4-2)
    {\small{ \sf Theorem \ref{the:4}: A symmetric independent MF sol. with decentralized information structure is near-optimal.}};
    \node[black,text width=3cm] at (-4.3,-5.5) (the5)
    {\small{\sf Theorems \ref{the:veri}--\ref{the:Approx-MK}: A finite MDP approximation for the representative DM MDP problem.}};
        \node[rectangle,purple,draw,fit=(the5),inner sep=0] {};
    \node[black,text width=3.4cm] at (3.5,-0.4) (the1-2)
    {\small{\sf Theorem \ref{the:2}: An optimal sol. exists that is exchangeable with the centralized information structure.}};
    \node [rectangle,purple,draw,fit=(the4-2),inner sep=0] {};
    \node [rectangle,purple,draw,fit=(the1-2),inner sep=0] {};
     \node[text width=8cm] at (8,1.5) 
    {\bf Infinite Horizon Discounted Teams};
    \fill[black!15] (7.6,0) ellipse (2.3 and 0.8);
    \fill[black!15] (7.6,-5) ellipse (2 and 0.7);
    \fill[black!15] (5.8,-2.5) ellipse (1.8 and 0.7);
     \node[text width=4cm] at (7.8,0.5) 
    {$N$-DMs Teams with};
    \node[text width=3cm] at (8.1,0.1) 
    {Centralized};
    \node[text width=4cm] at (7.9,-0.3) 
    {information structure};
    \node[text width=4cm] at (8.2,-5) 
    {Mean-field limit};
    \node[text width=3cm] at (5.8,-2.3) 
    {Measure-valued};
    \node[text width=3cm] at (6.7,-2.8) 
    {MDP};
    \draw[black, thick, ->]  (3,-2.5) -- (4,-2.5);
    \draw[black, thick, ->]  (8,-4.5) -- (8,-0.5);
    \draw[black, thick, <->]  (6.5,-0.6) -- (6,-1.9);
    \draw[black, thick, <-]  (6.5,-4.5) -- (6,-3.1);
    \node[black,text width=0.1cm] at (5.6,-3.4) 
    {\small{$N$}};
    \node[black,text width=0.1cm] at (5.9,-4.3) 
    {\small{$\infty$}};
    \draw[black, thick, ->]  (5.8,-3.6) -- (6,-4.1);

    \node[black,text width=2cm] at (3.8,-2.2) 
    {\small{$T\to\infty$}};
    \draw[blue, thick, ->]  (4.8,-5) -- (5.8,-5);
    \draw[blue, thick, ->]  (5,-0.4) -- (5.7,-0.2);
    \draw[blue, thick, <-]  (1.1,-5) -- (2.1,-5);
    \draw[blue, thick, <-]  (1,-0.1) -- (1.8,-0.2);
    \draw[blue, thick, ->]  (-3,-0.8) -- (-2.5,-0.5);
    \draw[blue, thick, ->]  (9,-0.8) -- (8.5,-0.5);
    \draw[blue, thick, ->]  (1.9,-1.5) -- (1.6,-2.2);
    \draw[blue, thick, ->]  (4.3,-1.4) -- (4.6,-2.1);
\node[black,text width=3cm] at (11.4,-5.5) (the7)
    {\small{ \sf Theorems \ref{the:veri}--\ref{the:Approx-MK}: A finite MDP approximation for the representative DM MDP problem.}};
    \node[rectangle,purple,draw,fit=(the7),inner sep=0] {};
\end{tikzpicture}
}
\caption{Main results for finite and infinite population teams.}
\label{fig:1}
\end{figure*}

The main contributions of the paper are summarized in Fig. \ref{fig:1}.

We first consider \PNT\ under randomized policies with the centralized information structure. We define various sets of randomized policies and their corresponding strategic measures on joint probability measures on states, actions of all DMs, and the empirical process on states. We first use the equivalent measure-valued MDP approach in \cite{bauerle2023mean,Ali-2023} to derive the value iteration process for \PNT\ with the centralized information structure. Following this MDP formulation, we show that an optimal Markov policy satisfying value iterations (see \eqref{eq:value-iter-T}) leads to an exchangeable joint distribution of states and actions among DMs induced by exchangeable Markov policies. This completely characterizes an optimal policy and establishes the existence of an optimal policy that is exchangeable and Markov for \PNT\ (see Theorem \ref{the:2}).  We emphasize that for this class of policies, randomization among DMs can still be correlated and asymmetric. For \PN, using the classical fixed point argument and characterizing the limit of \PNT\ as $T\to \infty$, we show that an optimal policy exists that is exchangeable and stationary. We emphasize that exchangeability does not imply symmetry, and hence, we cannot guarantee that an optimal policy is symmetric, conditionally independent, and of decentralized MF sharing type (that is, it only uses local states and actions, and the empirical measure of states; see \eqref{eq:IS-DEC}). 

Our main result for \PIT\ is to show that an optimal policy exists that is symmetric and independent among DMs, and of decentralized MF sharing type. Our approach is to characterize the limit of a sequence of optimal exchangeable policies for \PNT\ as $N\to \infty$. We first use a finite de Finetti type theorem by Diaconis and Friedman (see Theorem \ref{the:dia}), to construct a sequence of infinitely exchangeable probability measures on states and actions that are close in the total variation distance to those induced by optimal exchangeable policies of \PNT. We then use the compactness of the set of infinitely exchangeable distributions under the weak convergence topology to show the existence of a converging subsequence under the weak convergence topology for any finite marginals. Then, we utilize a theorem by Aldous (see Theorem \ref{the:ald}) and a finite de Finetti theorem by Diaconis and Friedman (see Theorem \ref{the:dia}), to show that a subsequence of optimal exchangeable policies for \PNT\ converges as $N\to \infty$ to a policy that is infinitely-exchangeable. Additionally, utilizing Theorem \ref{the:ald}, we show that the empirical measure on states and state-action pairs induced by optimal exchangeable policies of the finite population setting converges to the directing measures of states and state-action pairs induced by an infinitely exchangeable policy. By de Finetti's Theorem, this infinitely exchangeable policy is symmetric and independent, conditioned on common randomness (its directing measure), which is characterized by the limit in the distribution of the empirical state and action process induced by an optimal exchangeable policy of \PNT. We emphasize that this limiting policy is conditionally symmetric and independent and only utilizes the decentralized MF-sharing type information structure. We finally use value iterations for \PNT\ to show that the limiting policy is optimal for \PIT\ (see Theorem \ref{the:3}) and then establish the optimality of the representative McKean-Vlasov control problem. We establish similar results for \PIN\ using the classical fixed point argument and characterizing the limit of \PIT\ as $T\to \infty$. As a corollary to Theorem \ref{the:3}, we establish optimality of the McKean-Vlasov control problem in Corollary \ref{Cor:MV} and then provide a verification theorem in Theorem \ref{the:veri}.

Recall that for \PNT, we only establish the existence of an optimal policy that is exchangeable under the centralized information structure. In Theorem \ref{the:4}, we show near-optimality of an optimal symmetric and independent policy of the McKean-Vlasov control problem that utilizes only decentralized MF sharing type information structure for \PNT\ with a finite but large number of DMs. Finally, in Theorem  \ref{the:Approx-MK}, we provide a finite approximation model for the MDP of the representative DM by quantizing the measure-valued MDP for the McKean-Vlasov control problem.

\section{Sets of Randomized Policies.}\label{sec:Randomization}
In the following, we define various sets of randomized policies for the finite population as well as infinite population teams. We start with the finite population teams, and then we extend the definitions for the infinite population teams.

\subsection{Randomized policies for the finite population teams.}
We start by defining the set of jointly randomized policies with the centralized information structure in \eqref{eq:IS} for the finite population teams. 

A (jointly) randomized policy $\{{\pi}_t^N\}_{t\geq0}$ is a sequence of stochastic kernels ${\pi}_t^N$ with $${\pi}_t^N \in  \Pc\left(\prod_{i=1}^{N}\U\bigg|\prod_{k=0}^{t-1}\left(\prod_{i=1}^{N}(\X \times \U) \times \mathcal{P}(\X)\right)\times \prod_{i=1}^{N}\X \times \mathcal{P}(\X)\right),$$
     given history of centralized information  ${h}_{t}:=(x_{0:t}^{1:N},u^{1:N}_{0:t-1},\mu_{0:t}^N)$ for all $t\geq 0$. We denote the set of all these randomized policies by $\PPN$. We note that stochastic kernels $\{{\pi}_t^N\}_{t\geq0}$ are selected without imposing any constraint on randomization devices among DMs, that is, given the centralized information available to each DM, the randomization in actions might be correlated. 
    Furthermore, a policy $\{{\pi}_t^N\}_{t\geq0}$ is \emph{Markov} if 
    for any $t\geq 0$,
\begin{align*}
    \pi_{t}^N(u_{t}^{1:N}\in \cdot|h_t)=\pi_{t}^N(u_{t}^{1:N}\in \cdot|x_t^{1:N},\mu_t^N).
\end{align*}
In other words, the history is irrelevant given the current states and the empirical measure. 
We denote Markov policies within randomized policies by adding a sub-index $M$, as $\PM$.
A policy $\{{\pi}_t^N\}$ is \emph{stationary Markov} (or simply stationary) if it is Markov and time-invariant. We denote stationary policies within randomized policies by adding a sub-index $S$, as $\PS$. 

In the following, we introduce subsets of the randomized policies above by restricting how DMs select their randomization devices.

We first define randomized policies that are exchangeable. In the following, we recall the definition of \emph{exchangeability} for random variables.
\begin{definition}\label{def:exc}
Random variables $q^{1:m}$ defined on a common probability space are $m$-\emph{exchangeable} if for any permutation $\sigma$ of the set $\{1,\dots,m\}$ (denoted by $\sigma\in S_{m}$), $$\mathcal{L}\left(q^{\sigma(1):\sigma(m)}\right)=\mathcal{L}\left(q^{1:m}\right),$$ where $q^{\sigma(1):\sigma(m)}:=(q^{\sigma(1)},\ldots,q^{\sigma(m)})$ and $\mathcal{L}(\cdot)$  denotes the (joint) law of random variables. Random variables $(q^{1},q^{2},\dots)$ are \emph{infinitely-exchangeable} if finite distributions of $(q^{1},q^{2},\dots)$ and $(q^{\sigma(1)},q^{\sigma(2)},\dots)$ are identical for any finite permutation (affecting only finitely many elements) of $\mathbb{N}$. 
Additionally, stochastic kernel $\mathcal{Q}\in \mathcal{P}(\prod_{i=1}^{N} \mathbb{Y}|\prod_{i=1}^{N} \mathbb{X})$ is {permutation invariant} if
\begin{align*}
    \mathcal{Q}(y^{1:N}\in \cdot|x^{1:N})= \mathcal{Q}(y^{\sigma(1):\sigma(N)}\in \cdot|x^{\sigma(1):\sigma(N)})
\end{align*}
for all permutation $\sigma$ of $\{1, \ldots, N\}$.
\end{definition}

An exchangeable randomized policy $\{{\pi}_t^N\}_{t\geq 0}$ is a sequence of randomized policies ${\pi}_t^N$ such that
 \begin{align*}
\pi_{t}^N(u_{t}^{\sigma(1):\sigma(N)}\in \cdot|h_t^{\sigma})= \pi_{t}^N(u_{t}^{1:N}\in \cdot|h_t),
 \end{align*}
 for any permutation $\sigma$ of $\{1, \ldots, N\}$, where
 $h_{t}^{\sigma}:=(x_{0:t}^{\sigma(1):\sigma(N)},u_{0:t-1}^{\sigma(1):\sigma(N)},\mu_{0:t}^N)$ for all $t\geq 0$. We denote the set of all exchangeable policies by $\PEXN$. We note that the information structure utilized above remains centralized since DMs have access to the entire set of random variables in the history of information $h_{t}$ at any $t\geq 0$. Furthermore, the randomization devices among DMs can still be correlated at each time, but they need to be exchangeable (that is, permutation invariant). We can similarly define Markov and stationary policies within the exchangeable ones, denoted respectively by $\PEXNM$ and $\PEXNS$.

Next, we define another subset of randomized policies. We are interested in the case that DMs act independently (with private randomness and possibly common randomness), given a decentralized MF sharing information structure 
\begin{align}\label{eq:IS-DEC}
    I_{t, {\sf DEC}}^{i}:=\left\{x_{0:t}^{i},u_{0:t-1}^{i}, \mu_{0:t}^N\right\} \quad \forall t\geq1
\end{align}
with $I_{0, {\sf DEC}}^{i}:=\left\{x_{0}^{i}, \mu_{0}^N\right\}$.
Under this information structure, DMs have access only to their history of private states and actions, and also the history of the state's empirical process.
Due to the exchangeability of the team problem and the fact that there can be many exchangeable DMs within the team problem, we are interested in \emph{symmetry}, in the sense that every DM selects an identical (randomized) policy (possibly conditioned on common randomness). To this end, we define symmetric and independent policies with common randomness under the decentralized MF sharing information structure in \eqref{eq:IS-DEC}. A randomized policy  $\{{\pi}_t^N\}_{t\geq0}$ is symmetric and independent with common randomness if
    \begin{align*}
     \pi_{t}^N(u_{t}^{1:N}\in \cdot|h_t)=\int_{0}^{1} \prod_{i=1}^{N}\overline{\pi}_{t}^N(u_{t}^{i}\in \cdot|h_t^i,z)\nu_t(dz),
 \end{align*}
for $\nu_t\in \Pc([0,1])$ is the distribution of common randomness (independent from the intrinsic exogenous
system random variables) and for some stochastic kernel
$$\overline{\pi}_t^N\in \mathcal{P}\left(\U\bigg| \prod_{k=0}^{t-1}\left(\X\times \U\times \mathcal{P}(\X)\right)\times \X \times \mathcal{P}(\X) \times [0,1]\right),$$ where 
$h_{t}^i:=(x_{0:t}^{i},u^i_{0:t-1},\mu_{0:t}^N)$ for all $t\geq 0$ correspond to the private history of DM$^i$ at time $t$. We denote this set of policies by $\PCONS$. We note that in the above, randomized policies $\pi_{t}^N$ are identical among DMs, independent conditioned on the common randomness $z$, and decentralized (i.e., it utilizes the decentralized MF sharing information in \eqref{eq:IS-DEC}).  

We next define symmetric and independent policies. 
A randomized policy  $\{{\pi}_t^N\}_{t\geq0}$ is symmetric and independent if
\begin{align*}
     \pi_{t}^N(u_{t}^{1:N}\in \cdot|h_t)= \prod_{i=1}^{N}\overline{\pi}_{t}^N(u_{t}^{i}\in \cdot|h_t^i),
 \end{align*}
for some stochastic kernel
$$\overline{\pi}_t^N\in \mathcal{P}\left(\U\bigg| \prod_{k=0}^{t-1}\left(\X\times \U\times \mathcal{P}(\X)\right) \times \X\times \mathcal{P}(\X)\right).$$ We denote this set of policies by $\PPRNS$.

Symmetry has an important implication that allows us to propagate $(x_t^i, \mu_t^N)$  via \eqref{eq:dyn} only utilizing private history $h_{t}^i$. As we will show in the following lemma, the reason behind this is that under (conditionally independent) symmetric policies, pairs $(x_{t}^{i},\mu_t^N)$ are conditionally independent of $x_{t}^{1:N}$ given $x_{t-1}^i, u_{t-1}^{i},\mu_{t-1}^N$. This allows us to propagate $\mu_{t}^N$ over time using only $x_{t-1}^i, u_{t-1}^{i},\mu_{t-1}^N$ given symmetric policies $\overline{\pi}_{t}^N$. A similar result is established in \cite[Theorem 8]{yongacoglu2022independent} and \cite{Mahajan-CDC14} for problems with finite spaces of states and actions. We, however, provide an alternative proof technique by exploiting exchangeability and symmetry in policies. The proof of the following lemma is included in the Appendix \ref{APP:Lem-sym}.

\begin{lemma}\label{lem:sym}
    Consider dynamics in \eqref{eq:dyn}. Under symmetric randomized policies $\pmb{\pi}=\{{\pi}_{t}\}_{t\geq 0}\in \PCONS$, 
\begin{align}\label{eq:lem-eq1}
&Pr^{\pmb{\pi}}\left((x^{i}_{t+1},\mu_{t+1}^N)\in\cdot|x_{0:t}^{1:N},u_{0:t}^{1:N},z\right) = Pr^{\pmb{\pi}}\left((x^{i}_{t+1},\mu_{t+1}^N)\in\cdot|x_{t}^i, u_{t}^{i},\mu_{t}^N,z\right),
\end{align}
and under symmetric randomized policies $\pmb{\pi}=\{{\pi}_{t}\}_{t\geq 0}\in\PPRNS$
\begin{align}\label{eq:lem-eq2}
&Pr^{\pmb{\pi}}\left((x^{i}_{t+1},\mu_{t+1}^{N})\in\cdot|x_{0:t}^{1:N},u_{0:t}^{1:N}\right)= Pr^{\pmb{\pi}}\left((x^{i}_{t+1},\mu_{t+1}^{N})\in \cdot|x_{t}^i, u_{t}^{i},\mu_{t}^N\right).
\end{align}
\end{lemma}

The above lemma states that $\{x_{t}^{i}, \mu_{t}^{N}, u_{t}^{i}\}_{t\geq 0}$ is a controlled Markov chain under any symmetric policy. The main idea behind the proof of the above lemma is to show that under symmetric policies, for all $t\geq 0$, all states with the same empirical measure $\mu_{t}^N$ lead to possible actions with the joint state-action empirical measure of $\theta_t^N$. This implies that DM$^{i}$ can deduce the joint state-action empirical measure $\theta_t^N$ under the decentralized information structure, thanks to symmetry in policies. Then, we show that the joint conditional probability of $\mu_{t+1}^N$ and $x_{t+1}^{i}$ given $\mu_t^N,\theta_t^N, x_t^i$ and $\overline{\pi}_t^N$ remains the same as that given $x_{0:t}^{1:N}$ and $u_{0:t}^{N}$ with joint state-action empirical measures $\theta_t^N$ and marginals on states as $\mu_t^N$. This is crucial for our analysis as it allows us to propagate $x_{t+1}^{i}$ and $\mu_{t+1}^N$ under symmetric policies using only the decentralized MF sharing information structure.  

\subsection{Randomized policies for infinite population teams.}

We start by defining the set of jointly randomized policies with a centralized information structure for infinite population teams.

A (jointly) randomized policy $\{{\pi}_t\}_{t\geq0}$ is a sequence of stochastic kernels ${\pi}_t$ given history of the centralized information $\widehat{h}_{t}$ such that 
\begin{align}\label{eq:Pi-policy}
{\pi}_t \in  \Pc\left(\prod_{i=1}^{\infty}\U\bigg|\left(\prod_{k=0}^{t-1}\prod_{i=1}^{\infty}(\X\times \U) \right)\times \prod_{i=1}^{\infty}\X\right),
\end{align}
     where $\widehat{h}_{t}:=(\underline{x}_{0:t},\underline{u}_{0:t-1}):=(x_{0:t}^{1},u_{0:t-1}^{1},x_{0:t}^{2},u_{0:t-1}^{2}, \ldots)$ for all $t\geq 0$. We drop the sub-index $N$ for randomized policies for the infinite population teams. We denote this set of policies by $\Pi$. Furthermore, a policy $\{\pi_t\}_{t\geq0}$ is Markov if 
    for any $t\geq 0$,
\begin{align*}
    \pi_{t}(\underline{u}_{t}\in \cdot|\widehat{h}_t)=\pi_{t}(\underline{u}_{t}\in \cdot|\underline{x}_t).
\end{align*}
A policy $\{\pi_t\}_{t\geq0}$ is stationary Markov if it is Markov and time-invariant. We denote these sets of policies by $\Pi_{\sf M}$ and $\Pi_{\sf S}$, respectively.

Since our goal is to establish exchangeability or symmetry for an optimal policy of the infinite population MF teams, we are interested in the extension of the finite exchangeable policies to the infinite population limit. In particular, we want to define randomized policies that are infinitely exchangeable (see Definition \ref{def:exc}).  An infinitely-exchangeable randomized policy $\{{\pi}_t\}_{t\geq 0}$ is a sequence of randomized policies $\{{\pi}_t\}_{t\geq 0}$ such that
 \begin{align*}
      \pi_{t}(\underline{u_{t}^{\sigma}}\in \cdot|\widehat{h}_t^{\sigma})= \pi_{t}(\underline{u_{t}}\in \cdot|\widehat{h}_t), \forall t\geq 0
 \end{align*}
 for any finite permutation $\sigma$ of $x^{1:N}_{0:t},u^{1:N}_{1:t}$ among DMs for all $N\geq 1$. We define this set of policies by $\PEX$, and its corresponding Markov and stationary exchangeable policies by $\PEXM$ and $\PEXS$, respectively.

Next, we define symmetric and independent policies with common randomness under the decentralized MF sharing information structure. A randomized policy  $\{{\pi}_t\}_{t\geq0}$ is symmetric and independent with common randomness if
    \begin{align}
     \pi_{t}(\underline{u}_{t}\in \cdot|\widehat{h}_t)=\int_{0}^{1} \prod_{i=1}^{\infty}\overline{\pi}_{t}(u_{t}^{i}\in \cdot|\widehat{h}_t^i,z)\eta(dz),
 \end{align}
for some stochastic kernel
$$\overline{\pi}_t\in \mathcal{P}\left(\U\bigg| \prod_{k=0}^{t-1}\left(\X\times\U\times \mathcal{P}(\X)\right)\times \X\times \mathcal{P}(\X) \times [0,1]\right),$$ where 
$\widehat{h}_{t}^i:=(x_{0:t}^{i},u_{0:t-1}^i,\mu_{0:t})$ for all $t\geq 1$ and $\widehat{h}_{0}^i:=(x_{0}^{i},\mu_0)$ correspond to the private history of DM$^i$. In the above,  $\mu_{t}$ is the conditional distribution of  ${x}^i_t$ given the common randomness $z$, i.e., ${\mu_t}=\mathcal{L}({x}^i_t|z)\in \Pc(\X|[0,1])$ for all $t\geq 0$. We denote the set of these policies by $\PCOS$. We note that ${\mu_t}$ are identical among DMs since policies are identical conditioned on $z$.

We next define symmetric and independent policies. 
A randomized policy  $\{{\pi}_t\}_{t\geq0}$ is symmetric and independent if
    \begin{align*}
     \pi_{t}(\underline{u}_{t}\in \cdot|\widehat{h}_t)= \prod_{i=1}^{\infty}\overline{\pi}_{t}(u_{t}^{i}\in \cdot|\widehat{h}_t^i),
 \end{align*}
for some stochastic kernel
$$\overline{\pi}_t\in \mathcal{P}\left(\U\bigg| \prod_{k=0}^{t-1}\left(\X\times\U\times \mathcal{P}(\X)\right)\times \X\times \mathcal{P}(\X)\right)$$ with  ${\mu_t}=\mathcal{L}({x}^i_t)\in \Pc(\X)$ for all $t\geq 0$.  We denote the set of these policies by $\PPRS$. We again note that ${\mu_t}$ are identical among DMs due to symmetry in policies and dynamics. We additionally note that although policies for the infinite population limit are well-defined under the centralized information structure, due to a possible lack of symmetry or exchangeability, we cannot define policies under the decentralized MF-sharing information structure without an ergodicity condition. This is because without symmetry or exchangeability, the limit of the sequence of empirical measures $\pmb{\mu^N}$ might be ill defined.

\section{Existence of Exchangeable and Symmetric Optimal Policies.}\label{sec:Existence}

In this section, we establish the existence and structural results for both finite population teams as well as infinite population ones. To this end, we need to first introduce various sets of strategic measures.

\subsection{Strategic measures for teams with finite and infinite population DMs.}\label{sec:strategic measures}

We now introduce sets of strategic measures for finite population teams \PNT\ and \PN\ and the infinite population teams \PIT\ and \PIN. We start with strategic measures for \PN (those for \PNT\ can be defined similarly). 

We start by introducing the set of all strategic measures, denoted by $\LN$. A probability measure $${P_{\pi}^N}\in  \Pc\left(\prod_{t=0}^{\infty}\left(\prod_{i=1}^{N}\left(\X \times \U\right) \times \EMX\right)\right)$$ is a strategic measure induced by the sequence of randomized policies $\{{\pi}_t^N\}_{t\geq0}\in \PPN$ with $x_{0}^{i}\sim P_{0}$ for all $i=1,\ldots, N$ if for every $t\geq 0$:
\begin{align}
    &\int g({h}_{t-1},x^{1:N}_t, \mu_t^N){P_{\pi}^N}(d{h}_{t-1},dx^{1:N}_t,d\mu_t^N)\nonumber\\
    &=\int {P_{\pi}^N}(d{h}_{t-1})\int_{\prod_{i=1}^{N}\X} g({h}_{t-1},x^{1:N}_t,\mu_t^N)Pr(d\mu_t^N|x^{1:N}_t)  \prod_{i=1}^{N}{\mathcal{T}}(dx^{i}_t|x_{t-1}^{i},u_{t-1}^{i},\mu_{t-1}^N)\label{eq:prop-dyn}
\end{align}
for any continuous and bounded function $g$. Recall that ${h}_{t}:=(x_{0:t}^{1:N},u_{0:t-1}^{1:N},\mu_{0:t}^N)$ for all $t\geq 1$ and ${h}_{0}:=(x_{0}^{1:N},\mu_0^N)$. {In the above, we utilized the notation $P_{\pi}^{N}$ to denote the joint probability measure and its marginals on the corresponding random variables, e.g., $P_{\pi}^{N}(d{h}_{t-1})$ corresponds to the marginals of the joint probability measure $P_{\pi}^{N}$ on random variables in ${h}_{t-1}$.} 

We denote this set of randomized policies (strategic measures) by $\LN$ endowed with the weak convergence topology. A strategic measure 
${P_{\pi}^N}$ is a joint distribution on states, actions, and empirical measure processes $(\pmb{x^{1:N}},\pmb{u^{1:N}}, \pmb{\mu^N})$ that describes the evolution of states and its empirical measure given the sequence of randomized policies $\{{\pi}_t^N\}_{t\geq0}\in \PPN$. We emphasize that the information structure used for policies $\{{\pi}_t^N\}_{t\geq0}\in \PPN$ above is centralized.

We now introduce a subset of $\LN$ including only strategic measures ${P_{\pi}^N}$ induced by exchangeable randomized policies $\{{\pi}_t^N\}_{t\geq0}\in \PEXN$. Denote the set of exchangeable strategic measures by $\LEXN$ that is given by 
\begin{align*}
    \LEXN:=&\bigg\{{P_{\pi}^N}\in \LN\bigg| {P_{\pi}^N}\left(\left(\pmb{x}^{1:N},\pmb{u}^{1:N},\pmb{\mu^N}\right)\in \cdot\right)={P_{\pi}^N}\left(\left(\pmb{x}^{\sigma(1):\sigma(N)},\pmb{u}^{\sigma(1):\sigma(N)},\pmb{\mu^N}\right)\in \cdot\right) \quad\forall \sigma\in S_N \bigg\},
\end{align*}
where $S_N$ is the set of all permutations $\sigma$ of the set $\{1, \ldots, N\}$. This is because under exchangeable randomized policies $\{{\pi}_t^N\}_{t\geq0}\in \PEXN$, the states and actions are exchangeable due to the symmetry of dynamics in \eqref{eq:dyn}. Note that $\pmb{\mu^N}$ is permutation invariant independent of policies, i.e, under any permutation $(\pmb{x}^{\sigma(1):\sigma(N)},\pmb{u}^{\sigma(1):\sigma(N)})$, the empirical measure $\pmb{\mu^N}$ remains the same.

Next, we define the set of symmetric (decentralized) strategic measures with private and common randomness as follows:
\begin{align*}
    \LCONS:=&\bigg\{{P_{\pi}^N}\in \LN\bigg| {P_{\pi}^N}\left((\pmb{{x}^{1:N}},\pmb{{u}^{1:N}},\pmb{\mu^N})\in \cdot\right)=\int_{\cdot \times [0,T]} \prod_{i=1}^{N}\widetilde{{P_{\pi}^N}}(d\pmb{x^i},d\pmb{u^i},d\pmb{\mu^N}|z)\nu(dz)\bigg\},
\end{align*}
where $\nu\in \Pc([0,1])$ is the distribution of common randomness (independent from the intrinsic exogenous
system random variables) and $\widetilde{{P_{\pi}^N}}$ is a strategic measure in $$\Pc\left(\prod_{t=0}^{\infty}(\X \times \U\times \mathcal{P}(\X))\right)$$ induced by symmetric (identical) conditionally randomized policies $\{{\pi}_t^N\}_{t\geq0}\in \PCONS$ with a common randomness, i.e.,  for every $t\geq 0$
\begin{align}\label{eq:prop-x-mu}
    &\int g(h_{t-1}^i,x^{i}_t, \mu_t^N){\widetilde{P_{\pi}^N}}(dh_{t-1}^i,dx^{i}_t,d\mu_t^N|z)\nonumber\\
    &=\int {\widetilde{P_{\pi}^N}}(dh_{t-1}^i)\int_{\X} g(h_{t-1}^i,x^{i}_t,\mu_t^N)Pr(dx^{i}_t,\mu_t^N|x_{t-1}^i, u_{t-1}^{i},\mu_{t-1}^N,z)
\end{align}
for every continuous and bounded function $g$, where we recall that $h_{t}^i:=(x_{0:t}^{i},u_{0:t-1}^{i},\mu_{0:t}^N)$ for all $t\geq 0$ correspond to the private history of DM$^i$, and for each $t\geq 0$ actions $u_{t}^{i}$ are induced by identical policies $\overline{\pi}_{t}^{N}$ conditioned on $z$ given $h_{t}^i$. We note that in the above set of strategic measures $\overline{\pi}_{t}$ is identical among DMs, independent conditioned on the common randomness $z$, and decentralized (i.e., it utilizes the decentralized MF sharing information set in \eqref{eq:IS-DEC}). Since randomized policies are identical and independent conditioned on $z$, the state-action pairs $(x_{t+1}^{i},u_{t+1}^{i})$ are i.i.d. conditioned on $\mu_{t}^N$ and $z$. We note that in \eqref{eq:prop-x-mu}, since policies are symmetric $(x_t^i, \mu_t^N)$ can propagate via \eqref{eq:dyn} only utilizing private history $h_{t}^i$ thanks to Lemma \ref{lem:sym}.

We next define the set of all symmetric (identical) randomized policies induced by private (independent) randomness by $\LPRNS$ as follows 
\begin{align*}
    &\LPRNS:=\bigg\{{P_{\pi}^N}\in \LN\bigg|{P_{\pi}^N}\left((\pmb{{x}^{1:N}},\pmb{{u}^{1:N}},\pmb{\mu^N})\in \cdot\right)=\int_{\cdot} \prod_{i=1}^{N}\widetilde{{P_{\pi}^N}}(d\pmb{x^i},d\pmb{u^i},d\pmb{\mu^N})\bigg\}
\end{align*}
for some strategic measure $\widetilde{{P_{\pi}^N}}\subseteq  \Pc(\prod_{t=0}^{\infty}\X \times \U\times \EMX)$ induced by policies $\{\pi_{t}^{N}\}_{t\geq 0}\in \PPRNS$ which is defined as \eqref{eq:prop-x-mu} without conditioning on $z$. Using an analogous reasoning as the one used in Lemma \ref{lem:sym} for $\LCONS$, we can propagate $\mu_{t+1}^N$ over time using only $x_{t}^{i}, \mu_t^N$ and policies $\pi_{t}^N$ thanks to symmetry in policies.

We can similarly define strategic measures for the finite horizon problem \PNT. When it can lead to confusion, we denote strategic measures for \PNT\ by adding a sub-index $T$, i.e., ${P_{\pi,T}^N}$. We, however, do not include the sub-index for the sets of strategic measures (i.e., $L^{N}_{T}=L^N$ and so on).  

We now redefine the expected cost under the selection of strategic measures. For any ${P_{\pi}^N}\in \LN$, the discounted finite horizon and infinite horizon expected cost for finite population teams \PNT\ and \PN are given by
\begin{align*}
   J^N_T({P_{\pi,T}^N}):= \int \frac{1}{N}\sum_{i=1}^{N}\sum_{t=0}^{T-1}\beta^tc(x_{t}^{i},u_{t}^{i},\mu_{t}^N) {P_{\pi,T}^N}(d\pmb{x^{1:N}},d\pmb{u^{1:N}},d\pmb{\mu^N}),
\end{align*}
and 
\begin{align*}
  J^N({P_{\pi}^N}):=  \int \frac{1}{N}\sum_{i=1}^{N}\sum_{t=0}^{\infty}\beta^tc(x_{t}^{i},u_{t}^{i},\mu_{t}^N) {P_{\pi}^N}(d\pmb{x^{1:N}},d\pmb{u^{1:N}},d\pmb{\mu^N}).
\end{align*}

Next, we introduce various sets of strategic measures for the infinite population MF team \PIN\ (and \PIT). Since we aim to establish exchangeability or symmetry for an optimal policy of the infinite population MF teams, we are interested in extending the exchangeable strategic measures in $\LEXN$ to the infinite population limit. In particular, we want to define randomized policies that are infinitely exchangeable (see Definition \ref{def:exc}). Denote the set of infinitely exchangeable strategic measures $${P_{\pi}}\in \Pc\left(\prod_{t=0}^{\infty}\left(\prod_{i=1}^{\infty}\left(\X \times \U\right) \times \Pc(\X)\right)\right)$$ by $\LEX$ given by
\begin{align*}
    \LEX&:=\bigg\{{P_{\pi}}\in L\bigg| \forall N\in \mathbb{N}, {P_{\pi}}\left((\pmb{{x}^{1:N}},\pmb{{u}^{1:N}},\pmb{\mu})\in \cdot\right)={P_{\pi}}\left((\pmb{x}^{\sigma(1):\sigma(N)},\pmb{u}^{\sigma(1):\sigma(N)},\pmb{\mu})\in \cdot\right) \quad\forall \sigma\in S_N\bigg\},
\end{align*}
induced by randomized policies $\{\pi_{t}\}_{t\geq 0}\in \PEX$, where $\pmb{\mu}=\{\mu_{t}\}_{t\geq0}$ is the directing measure of $\underline{\pmb{x}}:=\{\underline{x}_t\}_{t\geq 0}$, i.e., the random variables $(x^1_t,x^2_t, \ldots)$ are identical and independently distributed (i.i.d.), conditioned on $\mu_t$ for all $t\geq 0$.  Before describing how the directing measures of $\mu_t$ propagate and also describing the information structure utilized in $\LEX$, we define the following set, which is identical thanks to a de Finetti representation theorem \cite[Theorem 1.1]{kallenberg2006probabilistic}, stating that any infinite exchangeable sequence of random variables is distributed as a mixture of i.i.d. sequences, and hence, the process of the directing measures $\pmb{\mu}$ in $\LEX$ propagates as the conditional law of private states. Also, the information used in $\LEX$ will be decentralized as in $\LCOS$. 

Next, we define the set of all symmetric (identical) randomized policies with
common and private randomness $\LCOS$ as follows
\begin{align*}
  \LCOS:=  \left\{{P_{\pi}}\bigg| {P_{\pi}}\left((\underline{\pmb{x}},\underline{\pmb{u}},\pmb{\mu})\in \cdot\right)=\int_{\cdot} \prod_{i=1}^{\infty}\widetilde{{P_{\pi}}}(d\pmb{x^i},d\pmb{u^i},d\pmb{\mu}|z)\nu(dz)\right\},
\end{align*}
induced by randomized policies $\{\pi_{t}\}_{t\geq 0}\in \PCOS$ for $\pmb{\mu}=\{\mu_t\}_{t\geq 0}$ with ${\mu_t}=\mathcal{L}({x}^i_t|z)\in \Pc(\X|[0,1])$ for all $t\geq 0$. We note that $\{\pi_{t}\}_{t\geq 0}$ are identical among DMs due to symmetry in policies and dynamics, where the dynamics are decoupled conditioned on $z$. 

We note that in contrast to \eqref{eq:prop-x-mu} where $\mu_t^N\in \EMX$ is the empirical measure of states at time $t$,  $\mu_t\in \Pc(\X|[0,1])$ is the law of state $x_{t}^{i}$ conditioned on the common randomness $z$. The law of $x_{t}^{i}$ is the same among all DMs due to symmetry. We note that the information structure used in $\LCOS$ is decentralized MF sharing conditioned on $z$ as in \eqref{eq:IS-DEC}.

We next define the set of all symmetric (identical) randomized policies with
only private randomness $\LPRS$ as follows
\begin{align*}
   \LPRS:= \left\{{P_{\pi}}\in L\bigg| {P_{\pi}}\left((\underline{\pmb{x}},\underline{\pmb{u}},\pmb{\mu})\in \cdot\right)=\int_{\cdot} \prod_{i=1}^{\infty}\widetilde{{P_{\pi}}}(d\pmb{x^i},d\pmb{u^i},d\pmb{\mu})\right\},
\end{align*}
induced by randomized policies $\{\pi_{t}\}_{t\geq0}\in \PPRS$ for $\pmb{\mu}=\{\mu_t\}_{t\geq 0}$ with ${\mu_t}=\mathcal{L}({x}^i_t)\in \Pc(\X)$ for all $t\geq 0$. Again due to symmetry $\{\mu_t\}_{t\geq 0}$ are identical among DMs. 

We can similarly define strategic measures for the finite horizon problem \PIT. We denote strategic measures for \PIT\ by adding a sub-index $T$, i.e.,   ${P_{\pi,T}}$. We similarly define the finite and infinite horizons discounted expected cost for any ${P_{\pi}}$ for \PIT\ and \PIN
\begin{align*}
J_T({P_{\pi,T}}):=\limsup_{N\to \infty}J_T^N({P_{\pi,T,N}})
\end{align*}
and 
\begin{align*}
J({P_{\pi}}):=\limsup_{N\to \infty}J^N({P_{\pi,N}})
\end{align*}
where ${P_{\pi,T,N}}$ and ${P_{\pi,N}}$ are the restrictions of ${P_{\pi,T}}$ and ${P_{\pi}}$, respectively to their first $N$ components.

\subsection{Existence of an exchangeable optimal policy for finite population teams.}

We first show the existence of an optimal policy for \PN\ and \PNT\ that is exchangeable. 

Following the approach in \cite{bauerle2023mean,Ali-2023}, we reformulate the problem as a measure-valued MDP. Let $\EM$ be the set of all empirical measures on states and actions $(x^{1:N}_t,u^{1:N}_t)$ for every $t\geq 0$ that are induced by any strategic measure ${{P_{\pi}^{N}}}\in L^N$. For any $\mu^N\in \EMX$, we define the new measure-valued action set by
\begin{align}\label{eq:act}
    U(\mu^N):=\left\{\theta^N\in \EM \mid \theta^N(\cdot \times \U)=\mu^N(\cdot)\right\}.
\end{align}

We define the new transition kernel
\begin{align}
    \mathbb{P}(\mu_{t+1}^N\in \cdot | \mu_{0:t}^N,\theta_{0:t}^N)&=\mathbb{P}\bigg(\mu_{t+1}^N\in \cdot \bigg|x_{t}^{1:N},u_{t}^{1:N}\:\text{s.t.}\:
     \frac{1}{N}\sum_{i=1}^{N}\delta_{(x_{t}^{i},u_{t}^{i})}=\theta_t^N, \frac{1}{N}\sum_{i=1}^{N}\delta_{x_{t}^{i}}=\mu_t^N\bigg)\label{eq:35}\\
    &:=\eta^N(\cdot|\mu_t^N,\theta_t^N)\label{eq:eta}
\end{align}
for a stochastic kernel $\eta^N\in \mathcal{P}\left(\mathcal{P}(\X)|\mathcal{P}(\X) \times \mathcal{P}(\X\times \U)\right)$, which describes a Markov dynamics that propagates the empirical measure process of the states $\mu_{t+1}^N$ given $\mu_t^N$ and $\theta_t^N$ for each $t\geq0$. We note that 
\begin{align}
    \mathbb{P}(\mu_{t+1}^N\in \cdot | x_{0:t}^{1:N},u_{0:t}^{1:N})
    &=\mathbb{P}(\mu_{t+1}^N\in \cdot | x_{t}^{1:N},u_{t}^{1:N})\nonumber\\
&=\mathbb{P}\bigg(\mu_{t+1}^N\in \cdot \bigg|x_{t}^{1:N},u_{t}^{1:N}\:\text{s.t.}\:\frac{1}{N}\sum_{i=1}^{N}\delta_{(x_{t}^{i},u_{t}^{i})}=\theta_t^N, \frac{1}{N}\sum_{i=1}^{N}\delta_{x_{t}^{i}}=\mu_t^N\bigg)\label{eq:38-T2}
\end{align}
where \eqref{eq:38-T2} follows from the following argument. By fixing $\theta_{t}^N$ and $\mu_{t}^N$, we can determine $x_{t}^{1:N}$ and $u_{t}^{1:N}$ up to their permutation. Due to exchangeability in the dynamics \eqref{eq:dyn}, we have
\begin{align*}
    \mathbb{P}(\mu_{t+1}^N\in \cdot | x_{t}^{1:N},u_{t}^{1:N})=\mathbb{P}(\mu_{t+1}^N\in \cdot | x_{t}^{\sigma(1):\sigma(N)},u_{t}^{\sigma(1):\sigma(N)})
\end{align*}
for any permutation $\sigma\in S_N$. This implies \eqref{eq:38-T2}.

For this new state and actions, the new running cost is given by
\begin{align*}
    \widetilde{c}(\mu^N,\theta^N)&:=\frac{1}{N}\sum_{i=1}^{N}c(x^i,u^{i}, \mu^N)=\int c(x,u,\mu^N)\theta^N(dx,du).
\end{align*}

Hence,
\begin{align*}
    \left(\EMX,\EM, \{U(\mu^N)|\mu^N\in \EMX\},\eta^N, \widetilde{c}\right)
\end{align*}
constitutes an MDP. By classical results for MDPs, we can restrict the search for an optimal policy to Markovian ones. As a result, we only define the space of Markov policies. Define the space of deterministic (Markov) policies at time $t\geq 0$ for the MDP formulation as a set of all measurable functions $g_t:\EMX\to \EM$ that maps the current measure-valued state $\mu_{t}^N$  to a measure-valued action $\theta_t^N$. We denote the set of these policies by $\mathbb{G}_t$. We can consider the finite-horizon discounted expected cost for any policy $\pmb{g_{T}}:=(g_0, \ldots, g_{T-1})\in \prod_{t=0}^{T-1}\mathbb{G}_t$ as follows
\begin{align*}
    J^{N}_{T}(\pmb{g_{T}})= \mathbb{E}^{\pmb{g_{T}}}\left[\sum_{t=0}^{T-1}\beta^t \widetilde{c}(\mu_t^N,\theta_t^N)\right].
\end{align*}
{Since the cost is non-negative and bounded, using the Fubini-Tonelli theorem, we can exchange the finite summation over $i$ with the infinite summation over $t$, i.e.,
\begin{align*}
    \frac{1}{N}\sum_{i=1}^{N}\sum_{t=0}^{\infty}\beta^tc(x_{t}^{i},u_{t}^{i},\mu_{t}^N)=\sum_{t=0}^{\infty}\beta^t\frac{1}{N}\sum_{i=1}^{N}c(x_{t}^{i},u_{t}^{i},\mu_{t}^N).
\end{align*}}
Hence, we define the infinite-horizon expected cost for any policy $\pmb{g}=(g_0, g_1, \ldots)\in \prod_{t=0}^{\infty}\mathbb{G}_t$ as 
\begin{align*}
    J^{N}(\pmb{g})= \mathbb{E}^{\pmb{g}}\left[\sum_{t=0}^{\infty}\beta^t \widetilde{c}(\mu_t^N,\theta_t^N)\right].
\end{align*}

For the finite-horizon expected cost, the value iteration equations  are given by
\begin{equation}
\left\{
\begin{aligned}
J^{N}_{T-1,T}(\mu^N)&:=\inf_{\theta^N\in U(\mu^N)}\widetilde{c}(\mu^N, \theta^N)\\
J^{N}_{t,T}(\mu^N)
&:=\inf_{\theta^N\in U(\mu^N)}\left\{\widetilde{c}(\mu^N, \theta^N)+\beta \int J^{N}_{t+1,T}({\mu}^{N}_{+}) \eta^N(d{\mu}^{N}_{+}\mid\mu^N, \theta^N)\right\}
\end{aligned}
\right.
\label{eq:value-iter-T}
\end{equation}
for all $t=0,\ldots T-2$ and every $\mu^N\in \EMX$. 

Using the monotone convergence theorem, the limit of $$J^N_{\infty}(\mu^N):=\lim\limits_{T\to \infty}J^{N}_{t,T}(\mu^N)$$ exists using a fixed-point argument for classical MDPs, and hence the value iteration equations for the infinite-horizon discounted expected cost is given by
\begin{align}
     J^{N}_{\infty}(\mu^N):=
     &\inf_{\theta^N\in U(\mu^N)}\left\{\widetilde{c}(\mu^N, \theta^N)+\beta \int J^{N}_{\infty}({\mu}^N_{+}) \eta^N(d{\mu}^N_{+}|\mu^N, \theta^N)\right\}\label{eq:value-iter-inf}
\end{align}
for every $\mu^N\in \EMX$.

Any admissible Markov policy $\pmb{g_T}\in \prod_{t=0}^{T-1}\mathbb{G}_t$ (or $\pmb{g}\in \prod_{t=0}^{\infty}\mathbb{G}_t$), induces a joint distribution on states, actions, and empirical measure of states that corresponds to a strategic measure ${P_{\pi}}\in L^N$. 

In the following theorem, we show that value iterations \eqref{eq:value-iter-T} (or \eqref{eq:value-iter-inf}) can, without any loss, be assumed to induce an optimal strategic measure that is exchangeable with Markovian randomized policies for  \PNT\ and stationary randomized policies for \PN. Proof of the following theorem is included in the Appendix \ref{APP-t2}.

\begin{theorem}\label{the:2}
    Consider \PNT\ and \PN. Let Assumption \ref{assump:1} hold. Then:
    \begin{itemize}
        \item [(i)] The value iterations \eqref{eq:value-iter-T}  admit a Markov solution $\pmb{g_{T}^{\star}}$ for \PNT\ that leads to an exchangeable joint distribution on $\pmb{x^{1:N}},\pmb{u^{1:N}},\pmb{\mu^N}$, induced by an exchangeable Markov randomized policy in $\PEXNM$. 
        \item [(ii)] The value iterations \eqref{eq:value-iter-inf}  admit a Markov stationary solution $\pmb{g^{\star}}$ for \PN\ that leads to an exchangeable joint distribution on $\pmb{x^{1:N}},\pmb{u^{1:N}},\pmb{\mu^N}$, induced by an exchangeable stationary randomized policy in $\PEXNS$. 
    \end{itemize}
\end{theorem}

Theorem \ref{the:2} establishes the existence of a globally optimal policy (strategic measure) ${P_{\pi}^{N\star}}$ that is exchangeable among all strategic measures in $L^N$, i.e., ${P_{\pi}^{N\star}}\in \LEXN$. Additionally, it provides further structural results for optimal randomized policies and shows that it is Markovian for \PNT\ and stationary for \PN. 

The idea of the proof is as follows: consider the measure-valued MDP formulation. We construct an exchangeable policy by averaging over all its permutations that perform the same as the given strategic measure. Specifically, we start at $T-1$, and we construct an exchangeable joint distribution on state and action pairs $(x_{T-1}^{1:N},u_{T-1}^{1:N})$ by averaging over all permutations of that induced by an optimal policy $g_{T-1}^{\star}$. We show that this exchangeable distribution attains the same team performance, and hence, it is optimal. Sequentially, we can show that this is true for every $t=0, \ldots, T-1$, and hence, we can conclude that an optimal strategy $\pmb{g^{\star}_{T}}$ induces an exchangeable joint distribution on state and action process $\pmb{x^{1:N}},\pmb{u^{1:N}}$ without any loss of optimality.  Weak-Feller continuity of the measure-valued MDP together with continuity of the cost then guarantees the existence of an optimal Markov policy $\pmb{g^{\star}_{T}}$, which completes the proof. 

In the following, we present a counterexample showing that  $\{g^{\star}_{t}\}_{t\geq 0}$ cannot be realized by any symmetric and independent policy under MF information sharing
$\pi^i_t(\cdot|x^i_t,\mu_t^N) = \pi^{\star}_t(\cdot|x^i_t,\mu_t^N),$ for $i=1,\cdots,N$ and $t\geq 0$. We note that a related example is provided in \cite[Section I.D]{Mahajan-CDC14} showing that symmetric deterministic policies with MF sharing information structure might not be optimal for the finite population teams. 

\begin{example}
Consider the following dynamics:
\[x^i_{t+1} = u^i_t, \qquad t\geq 0,\: i=1,2,\cdots,N.\]
We take $\mathbb{X}=\mathbb{U} = \{0,1\}$ (identical for each DM). Let the empirical measure be given with
\[\mu^N_t(z) = \frac{1}{N} \sum_{i=1}^N 1_{\{x^i_t=z\}}, \quad z \in \mathbb{X}.\]
Suppose that the cost to optimize is:
\begin{align}\label{costUnif}
E^{\pmb{\gamma^{1:N}}}\left[ \frac{1}{N}\sum_{t=0}^{T-1}\sum_{i=1}^{N} \sum_{z \in \mathbb{X}}\left|\mu^N_t(z) - \frac{1}{N}\right|^2\right].
\end{align}
where $\pmb{\gamma^{1:N}}$ is the set of all deterministic admissible policies. Here, we note that in (\ref{costUnif}) the goal is to make the empirical measure as uniform as possible.

Now, let us take $N=T=2$ and the initial joint state to be $x^1_0=x^2_0=1$. Then, under a symmetric (identical) policy which uses the MF information structure:
\[u^i_t=\gamma^i_t(x^i_t,\mu_t^N) = \gamma_t(x^i_t,\mu_t^N),\]
with $\mu_0^N=[0 \quad 1]$ and $x^1_0=x^2_0=1$ and fixed $\gamma_t$, we see that the MF term at time $t=1$ will always stay away from the uniform measure $[\frac{1}{2} \quad \frac{1}{2}]$ with arbitrarily high probability; if deterministic, it will always be a full mass on either $0$ or $1$, and if we allow for randomization, it will be a positive measure on all possible combinations of $(0,0), (0,1), (1,0), (1,1)$. The cost will be uniformly larger than $\frac{1}{2}$ as the second stage cost (at $t=1$) will be non-zero under any symmetric policy. If we allow for a small perturbation around identical policies, the cost will still be large, and if we allow for a large enough randomization, the probability of error will be bounded away from zero for the asymmetric empirical measures. In particular, the optimal symmetric randomized policy is to select the actions with probabilities $\frac{1}{2}$ leading to a cost of $0.5 +0.25=0.75$. 

On the other hand, if we allow for asymmetry in the policy 
\[\gamma^1_0(1,\mu_0^N)=1,\qquad \gamma^2_0(1,\mu_0^N)=0, \]
we can attain a cost of $0$ at time $1$, and the total cost will be $0.5$. 
\end{example}

We now provide an explanation for this counterexample. While the structural result on the control topology analysis leading to the action space
in \eqref{eq:act} is valid, the set of attainable actions under asymmetric policies is strictly larger than the set of attainable actions under symmetric policies. In the above, clearly, the joint measure obtained by
 \[P^{\gamma}(x^1_0=1,u^1_0=0)=1,\quad P^{\gamma}(x^2_0=1,u^2_0=0)=1\]
  is attainable under some policy map, and the map points to an optimal joint measure as an optimal action. However, this joint measure (optimal action) is not attainable under symmetric policies that use local state information and identical maps. 

\subsection{Existence of a symmetric decentralized optimal policy for infinite population teams.}

In the following theorem, we show the existence of an optimal policy for  \PIT\ and \PIN\ that is symmetric, independent, and decentralized. Our method uses a sequence of optimal exchangeable strategic measures for \PNT\ and \PN shown to exist in Theorem \ref{the:2}. The proof of Theorem \ref{the:3} is provided in the Appendix \ref{APP-t3}.

\begin{theorem}\label{the:3}
Consider  \PNT,\PN,\PIT, and \PIN. Let Assumption \ref{assump:1} hold. Then:
\begin{itemize}
    \item [(i)] A sequence of exchangeable optimal strategic measures $\{{P_{\pi,T}^{N\star}}\}_N\subseteq \LEXN$ for \PNT\ obtained via \eqref{eq:value-iter-T} ($\{{P_{\pi}^{N\star}}\}_N\subseteq \LEXN$ for \PN\  obtained via \eqref{eq:value-iter-inf}) converges through a subsequence to a symmetric and conditionally independent policy ${P_{\pi,T}^{\star}}\in \LCOS$ that is globally optimal for \PIT\ (${P_{\pi}^{\star}}\in \LCOS$ that is globally optimal for \PIN).
    \item [(ii)] There exists an optimal strategic measure ${P_{\pi,T}^{\star}}$ for \PIT\ (${P_{\pi}^{\star}}$ for \PIN) that is decentralized (with MF information sharing), symmetric, and independent (i.e., ${P_{\pi, T}^{\star}}\in \LPRS$ and ${P_{\pi}^{\star}}\in \LPRS$) with Markovian policies for \PIT\ and stationary policies for \PIN. 
\end{itemize}
\end{theorem}

We now shed some light on the proof of Theorem \ref{the:3}.  We first use a finite de Finetti type theorem by Diaconis and Friedman (see Theorem \ref{the:dia} in the Appendix), to construct a sequence of infinitely exchangeable probability measures on states and actions that are close in the total variation distance to those induced by optimal exchangeable strategic measures of \PNT.  Next, we utilize a de Finetti-type convergence theorem by Aldous (see Theorem \ref{the:ald} in the Appendix) and a theorem by Diaconis and Friedman (see Theorem \ref{the:dia} in the Appendix), to show that a subsequence of optimal exchangeable strategic measures of \PNT\ converges as $N\to \infty$ to a strategic measure in $\LEX$ (through a subsequence). By de Finetti's representation Theorem, this infinitely exchangeable strategic measure is symmetric and independent among DMs conditioned on common randomness, which is its directing measure. Furthermore, we show that its directing measure can be characterized by the limit in the distribution of the empirical state process induced by optimal exchangeable strategic measures of \PNT. This limiting strategic measure is conditionally symmetric and independent and only utilizes the decentralized MF sharing information structure in \eqref{eq:IS-DEC}. We finally use value iterations \eqref{eq:value-iter-T} for \PNT\ to show that the limiting strategic measure is optimal for \PIT. We establish similar results for \PIN\ using the classical fixed point argument and characterizing the limit of \PIT\ as $T\to \infty$. Part (ii) of Theorem \ref{the:3} follows from the fact that the directing measure of the limiting policy can be characterized by the weak limit of the empirical measures on joint states and actions.

We now present the following result as a corollary to Theorem \ref{the:3}.
\begin{corollary}\label{cor:1}
    Consider \PIT\ and \PIN. Let Assumption \ref{assump:1} hold. For \PIT\ (\PIN), an optimal solution under the decentralized MF-sharing information structure in \eqref{eq:IS-DEC} is optimal for \PIT\ (\PIN) within any policy with the centralized information structure in \eqref{eq:IS}.
\end{corollary}
\begin{proof}
    The proof follows from Theorem \ref{the:3} since we can restrict the search for an optimal strategic measure to $\LPRS$, which only utilizes the decentralized information $I_{t, {\sf DEC}}^{i}$ in \eqref{eq:IS-DEC} for $t\geq 0$.
\end{proof}

\section{Optimality of the McKean-Vlasov  Representative DM Optimal Control Problem,  Verification Theorem, and Value Iteration.}\label{sec:verif-0}
As an important implication of the results of the previous section and the proof of Theorem \ref{the:3}, we can rigorously justify that \PIT\ can be equivalently formulated as a single-agent representative DM problem. The following corollary to Theorem \ref{the:3}, establishes optimality of the representative McKean-Vlasov optimal control problem.     Proof of the following corollary is included in Appendix \ref{APP:CorMV}. 
      
\begin{corollary}\label{Cor:MV}
Consider the single-agent representative DM, with dynamics of the representative DM (with index $R$) given by 
\begin{align*}
x_{t+1}^{R}=f(x_t^{R},u_{t}^{R},\mu_{t}^{R},w_{t}^{R}),
\end{align*}
with $\mu^R_t=\mathcal{L}(x_{t}^{R})$, and the finite horizon cost
\begin{align*}
    J^{R}_{T}(P^R)=\int \sum_{t=0}^{T-1}\beta^tc(x_t^{R},u_{t}^{R},\mu_{t}^{R})P^R(d\pmb{x^{R}},d\pmb{u^{R}},d\pmb{\mu^R}),
\end{align*}
or the infinite horizon discounted cost 
\begin{align*}
    J^{R}(P^R)=\int \sum_{t=0}^{\infty}\beta^tc(x_t^{R},u_{t}^{R},\mu_{t}^{R})P^R(d\pmb{x^{R}},d\pmb{u^{R}},d\pmb{\mu^R})
\end{align*}
where $P^R$ is a strategic measure, joint distribution of path variables $\pmb{x^{R}},\pmb{u^{R}},\pmb{\mu^R}$, for the representative DM induced by a sequence of policies $\{\pi_{t}^R\}_{t\geq0}$ under the MF sharing information structure (which can use $x^{R}_{0:t},u_{0:t-1}^{R}, \mu_{0:t}^{R}$). 
Under Assumption \ref{assump:1}, any optimal solution of the above problem is optimal for \PIT\ (or \PIN). 
\end{corollary}

The above problem is a single-agent problem with a consistency condition $\mu^R_t=\mathcal{L}(x_{t}^{R})$--also known as \emph{McKean-Vlasov optimal control problem}. We can reformulate this as a measure-valued MDP where the dynamics with the transition kernel
\begin{align*}
    \mu_{t+1}^{R}\sim \eta(\cdot|\mu_{t}^{R},\theta_t^{R})
\end{align*}
with \begin{align*}
\eta(\cdot \mid \mu_t^{R},\theta_t^{R}):=P\!\left(\mu_{t+1}^{R} \in \cdot \mid \mu_t^{R},\theta_t^{R}\right)
&=
\delta_{\{F(\mu_t^{R},\theta_t^{R})\in \cdot\}}
\end{align*}
where
\begin{align*}
F(\mu_t^{R},\theta_t^{R})(\cdot)&:=
\int_{\cdot}
\mathcal{T}\!\left(dx^{R}_{t+1}\mid x_t^{R},u_t^{R},\mu_t^{R}\right)\,
\theta_t^{R}(dx_t^{R},du_t^{R}).
\end{align*}
Also, actions are $\theta^R_t=\mathcal{L}(x_{t}^{R},u_{t}^{R})$, and the running cost is given by
\begin{align*}
\widetilde{c}(\theta_t^{R},\mu_{t}^{R})=\int {c}(x_{t}^{R},u_{t}^{R},\mu_{t}^{R}) \theta_t^{R}(dx_{t}^{R},du_{t}^{R}).
\end{align*}

We can reformulate both finite and infinite horizon MDPs for the representative DM. Related to the above conclusion, \cite[Theorem 4.6]{bauerle2023mean} established that the limsup of the value function for the finite population MDP is equal to the value function of a limiting MDP as $N\to \infty$ using the compactness of state and action space and continuity of the MDP. However, the structure of an optimal policy has not been established. Here, however, we use exchangeability and de Finetti theorem to show the existence of a symmetric optimal solution for the infinite population problem and certify the optimality of value iterations for the McKean-Vlasov formulation.

We next provide the following verification theorem.
\begin{theorem}\label{the:veri}
    Consider the representative DM's MDP under Assumption \ref{assump:1}. 
    \begin{itemize}
\item [(i)] For the finite horizon, an optimal deterministic Markov policy $\pmb{g^{\star}_T}$ exists and satisfies the value iterations \begin{equation}
\left\{
\begin{aligned}
J_{T-1,T}(\mu)&=\inf_{\theta\in U(\mu)}\widetilde{c}(\mu, \theta)=\widetilde{c}(\mu, g^{\star}_{T-1}(\mu))\\
J_{t,T}(\mu)&=(\mathbb{T}(J_{t+1,T}))(\mu), \qquad \forall \mu\in \mathcal{P}(\X)
\end{aligned}
\right.
\label{eq:finite-horizon-dp}
\end{equation} for $t=0,\ldots, T-2$, where 
\begin{align*}
&\mathbb{T}(\nu)(z):=\inf_{\theta\in U(z)}\left\{\widetilde{c}(z, \theta)+\beta \int \nu({z}_{+}) \eta(d{z}_{+}|z, \theta)\right\}\nonumber. 
\end{align*}
Furthermore, an optimal policy $\pmb{g^{\star}_T}$ leads to a sequence of (possibly randomized) Markov policies $\{\pi_{t}^{\star}\}_{t\geq0}$ with $\pi_{t}^{\star}\in \mathcal{P}\left(\U\mid \X \right)$ for the representative DM under MF information sharing.
\item [(ii)] For the infinite horizon, an optimal deterministic stationary Markov policy $\pmb{g^{\star}}$ exists and satisfies the value iterations
\begin{equation*}
\left\{
\begin{aligned}
J_{\infty}(\mu)&=\widetilde{c}(\mu, g^{\star}(\mu))+\beta \int J_{\infty}({\mu}_{+})\, \eta(d{\mu}_{+}\mid\mu, g^{\star}(\mu))\\
J_{\infty}(\mu)&=(\mathbb{T}(J_{\infty}))(\mu), \qquad \forall \mu\in \mathcal{P}(\mathbb{X}).
\end{aligned}
\right.
\end{equation*}
Furthermore, an optimal policy $\pmb{g^{\star}}$ leads to a sequence of (possibly randomized) stationary Markov policies $(\pi^{\star},\pi^{\star},\ldots)$ for the representative DM with $\pi^{\star}\in \mathcal{P}\left(\U\mid \X \right)$ under MF information sharing.
\end{itemize}
\end{theorem}
\begin{proof}
Following an analogous argument as that in the proof of Theorem \ref{the:3} using the generalized dominated convergence theorem \cite{serfozo1982convergence}  (or \cite[Theorem 3.5]{langen1981convergence}), the transition kernel $\eta$ is weakly continuous, and $\widetilde{c}$ is continuous using Assumption \ref{assump:1}(ii)-(iii). Hence, the proof of the theorem follows from the measurable selection and verification theorem for classical MDPs (see e.g., \cite{HernandezLermaMCP}) since Assumption \ref{assump:1}(i) yields that its state space $\mathcal{P}(\X)$ is compact, and its action space 
$$U(\mu)=\left\{\theta_{\pi}\in \mathcal{P}(\mathbb{X} \times \mathbb{U})|\theta_{\pi}(\cdot) =\int\pi(du|x)\mu(dx)\right\}$$
is compact under the Young topology.
\end{proof}

\section{Approximations for Large Population Teams and the Representative DM Problem.}\label{sec:approximation}
In this section, we provide an approximation result for the large $N$-DM problem as well as finite approximations for the representative DM problem.

\subsection{Approximate symmetric decentralized optimal policies for large population teams.}

In Theorem \ref{the:2}, we showed that an optimal policy exists for \PNT\ (and \PN) that is exchangeable. This optimal policy, however, might not be symmetric and independent. In the following theorem, as an implication of our analysis in Theorem \ref{the:3}, we show that a symmetric and independent optimal solution of \PIT\ (and \PIN) with decentralized MF sharing information structure is near-optimal for \PNT\ (or \PN) with a large number of DMs among all policies under the centralized information structure. The proof of the following theorem is included in the Appendix \ref{APP:the4}.

\begin{theorem}\label{the:4}
    Consider \PNT\ and \PIT\ (\PN\ and \PIN). Let Assumption \ref{assump:1} hold. A symmetric independent optimal solution for \PIT\ (\PIN)  under the decentralized MF-sharing information structure in \eqref{eq:IS-DEC} is $\epsilon_{N}$-optimal for \PNT\ (\PN) within any policy with the centralized information structure in \eqref{eq:IS} with $\epsilon_{N}\to 0$ as $N\to \infty$.
\end{theorem}

The above theorem yields that a symmetric and independent solution of \PIT\ (or \PIN) with the decentralized MF sharing information structure is approximately optimal for \PNT\ (\PN) with a finite but large number of DMs. This allows us to use the representative DM formulation for \PIT\ (or \PIN) to obtain an approximate solution for  \PNT\ (\PN) when $N$ is large. 

Although we will not formally obtain a bound for the approximation error (i.e., $\epsilon_N$), we briefly discuss the rate of convergence in the following. Suppose that $P^{N\star}_{T}$ is an optimal strategic measure for \PNT. By Theorem \ref{the:2}, without any loss, we can assume that $P^{N\star}_{T}$ is exchangeable, i.e., $P^{N\star}_{T}\in \LEXN$. Our analysis in the proof of Theorems \ref{the:3} and \ref{the:4} shows that we can construct a sequence of strategic measures $\widetilde{P}^{N\star}_{T}\in \LEX=\LCOS$ such that their marginals to $N$ first components are close in the total variation to $P^{N\star}_{T}\in \LEXN$. We also showed that both $\widetilde{P}^{N\star}_{T}$ and  $P^{N\star}_{T}$ converge in distribution to $\widehat{P}^{\star}_{T}\in \LCOS$, an optimal strategic measure of \PIT\ with a common randomness as a weak limit of the sequence empirical measures of state-actions. Therefore, the limiting policy can be realized as a strategic measure in ${P}^{\star}_{T}\in\LPRS$ under MF information sharing. Also, the team's cost under $\widetilde{P}^{N\star}_{T}$ and  $P^{N\star}_{T}$ converge to that under $P^{\star}_{T}$. Therefore, the approximation error can be bounded as follows:
\begin{align*}
    \epsilon_{N}&=|J_{T}^{N}(P^{N\star}_{T})-J_{T}^{N}(P^{\star}_{T})|\leq |J_{T}^{N}(P^{N\star}_{T})-J_{T}^{N}(\widetilde{P}^{N\star}_{T})| + |J_{T}^{N}(\widetilde{P}^{N\star}_{T})-J_{T}^{N}(P^{\star}_{T})|.
\end{align*}
Since $\widetilde{P}^{N\star}_{T}$ and $P^{\star}_{T}$ are both symmetric and independent, a form of the law of large numbers can be used to obtain the rate of convergence for the empirical measure of states under (conditionally) i.i.d. randomized policies. This corresponds to a so-called \emph{propagation of chaos} in the literature, in which the rate is often of order $\mathcal{O}(1/\sqrt{N})$. This gives us an asymptotic bound for the second term above under some Lipschitz continuity conditions on the dynamics and cost. However, for bounding the first term, we need a rate of convergence for the empirical measure of states under exchangeable randomized policies, which requires further conditions and is left for future work. 

\subsection{Finite approximation for the representative DM problem and learning.}

We now provide an algorithm to obtain a near-optimal solution for the McKean-Vlasov control problem of the representative DM. We use a finite approximation approach for MDPs developed in \cite{saldi2018finite,saldi2017asymptotic} and also  \cite{kara2023q} to approximate the representative DM MDP. Finite approximation models for finite population MDPs have also been studied in \cite{Ali-2023} using the MDP connections in \cite{bauerle2023mean}, by quantizing the state space. In addition, considering symmetric policies for the infinite population MF problem, a finite approximation MDP is obtained for the limiting MDP \cite{Ali-2023} by discretizing the space of probability measures on states of the representative DM. Here, we focus on finite approximation of the MDP for the representative DM thanks to Theorem \ref{the:3} and Corollary \ref{Cor:MV}. We discretize the measure-valued state space as well as the actions equipped with the Young topology using the weak Feller property of the MDP and compactness of the underlying state and action spaces. 

First, we note that the above representative DM's MDP under Assumption \ref{assump:1} leads to a weak Feller kernel with the state space $\mathcal{P}(\X)$ and the action space 
$$U(\mu)=\left\{\theta_{\pi}\in \mathcal{P}(\mathbb{X} \times \mathbb{U})|\theta_{\pi}(\cdot) =\int\pi(du|x)\mu(dx)\right\},$$
which is compact under the Young topology. 

Since $(\mathbb{X},d_{\mathbb{X}})$ is a compact metric space, for any $x\in \mathbb{X}$, we can find a sequence $(\{x_{n,i}\}_{i=1}^{k_n})_{n\geq 1}$ of finite grids in $\mathbb{X}$ such that
\begin{align*}
    \min_{i=1, \ldots, k_n} d_{\mathbb{X}}(x,x_{n,i})<\frac{1}{n} \qquad \forall n\geq 1, 
\end{align*}
where $k_n$ is the size of the finite grid. We define the finite state space $\{x_{n,1},\ldots,x_{n,k_n}\}$ by $\widehat{\mathbb{X}}_n$, which is the quantization of $\mathbb{X}$ under nearest-neighbor quantizer $Q_n$ given by
\begin{align*}
Q_n(x)=\argminA_{x_{n,i}\in \widehat{\mathbb{X}}_n}\:d_{\mathbb{X}}(x,x_{n,i})\quad \forall x\in \mathbb{X}.
\end{align*} 
Let  $A_{n,i}=\{x\in \mathbb{X}|Q_n(x)=x_{n,i}\}$. This allows us to define the restriction of any $\mu\in \mathcal{P}(\mathbb{X})$ to $ A_{i,n}$ as $$\widehat{\mu}_{n,i}(\cdot)=\frac{\mu(\cdot)}{\mu(A_{n,i})}$$
provided that $\mu(A_{n,i})>0$ for all $i$ and $n$. Since $\X$ is compact, for any $\mu\in \mathcal{P}(\mathbb{X})$, there exists a sequence $\{\widehat{\mu}_n\}_n$ with $\widehat{\mu}_n\in \mathcal{P}(\widehat{\mathbb{X}}_n)$  that converges weakly to $\mu$ as $n\to \infty$. 

Similar to the state space, we quantize the action space and denote it by $\widehat{\mathbb{U}}_n$. Using this, we define 
\begin{align*}
\widehat{U}_n(\widehat{\mu})=\left\{\widehat{\theta}_{\pi,n}\in \mathcal{P}(\widehat{\mathbb{X}}_n \times \widehat{\mathbb{U}}_n)\bigg{|}\: \widehat{\theta}_{\pi,n}(\cdot)=\int\widehat{\pi}_{n}(du|x)\widehat{\mu}(dx)\right\}
\end{align*}
with $\widehat{\mu}\in\mathcal{P}(\widehat{\mathbb{X}}_n)$ and $\widehat{\pi}_n\in \mathcal{P}(\widehat{\mathbb{U}}_n|\widehat{\mathbb{X}}_n)$.
If $\mu$ is non-atomic, then $ \widehat{U}_n(\widehat{\mu})$ is dense in $U(\mu)$ (see e.g., \cite[Proof of Theorem 5.1]{yuksel2023borkar} \cite[Theorem 3.1]{beiglbock2018denseness}, \cite[Theorem 3]{milgrom1985distributional}, \cite{balder1997consequences}, \cite{borkar1988probabilistic}, or \cite[Chapter 7]{castaing2004young}). Hence, for any $\theta_{\pi}\in U(\widehat\mu)$, we can find a sequence $\{\widehat{\theta}_{\pi,n}\}_n$ (with  $\widehat{\theta}_{\pi,n}\in\widehat{U}_n(\widehat\mu)$) that converges to $\theta_{\pi}$ under the Young topology as $n\to \infty$. For any $\widehat{\mu}_n\in\ \mathcal{P}(\widehat{\mathbb{X}}_n)$ and  $\widehat{\theta}_{\pi,n}\in \widehat{U}_n(\widehat\mu)$, we define the transition kernel $\widehat{\eta}_n\in \mathcal{P}(\mathcal{P}(\widehat{\mathbb{X}}_n)|\mathcal{P}(\widehat{\mathbb{X}}_n)\times \mathcal{P}(\widehat{\mathbb{X}}_n \times \widehat{\mathbb{U}}_n))$ as
\begin{align*}
   \widehat{\eta}_n(\widehat{\mu}^{+}_{n,i}|\widehat{\mu}_n, \widehat{\theta}_{\pi,n}):=\eta(A_{i,n}|\widehat{\mu}_n, \widehat{\theta}_{\pi,n}) 
\end{align*}

We can also write the stage-wise cost as
\begin{align*}
   \widehat{c}(\widehat{P}_{\pi,n},\widehat{\mu}_{n})=\int c(x,u, \widehat{\mu}_n)  \widehat{\pi}_n(du|x)\widehat{\mu}_n(dx).
\end{align*}

As a result, 
\begin{align*}
 \left(\mathcal{P}(\widehat{\mathbb{X}}_n),\mathcal{P}(\widehat{\mathbb{X}}_n \times \widehat{\mathbb{U}}_n), \{\widehat{U}_n(\widehat{\mu})|\widehat{\mu}\in \mathcal{P}(\widehat{\mathbb{X}}_n)\},\widehat\eta_n, \widehat{c}\right)
\end{align*}
constitutes a finite MDP denoted by ${\sf MDP_n}$.

We can write the value iterations \begin{equation*}
\left\{
\begin{aligned}
\widehat{J}_{T-1,T}^n(\widehat{\mu})&=\inf_{\theta\in \widehat{U}_n(\widehat{\mu})}\widehat{c}(\widehat{\mu}, \theta)\\
\widehat{J}_{t,T}^n(\widehat{\mu})&=(\mathbb{T}_{n}(\widehat{J}_{t+1,T}^n))(\widehat{\mu}).
\end{aligned}
\right.
\end{equation*}
for a Bellman operator $\mathbb{T}_{n}$  on $\mathcal{P}(\widehat{\X}_n)$ given by
    \begin{align*}
&(\mathbb{T}_{n}(\widehat{\nu}))(\widehat{\mu}):=\inf_{\theta\in \widehat{U}_n(\widehat{\mu})}\left\{\widehat{c}(\widehat\mu, \theta)+\beta \int \widehat{\nu}({\mu}_{+}) \widehat\eta_n(d{\mu}_{+}|\widehat\mu, \theta)\right\}\nonumber.
\end{align*}

We can similarly define the value iteration for the infinite horizon discounted cost as follows: $$\widehat{J}_{\infty}^n(\widehat{\mu})=(\mathbb{T}_{n}(\widehat{J}_{\infty}^n))(\widehat{\mu}).$$ 

In the following theorem, we show that an optimal solution of ${\sf MDP_n}$ is approximately optimal for the MDP of the representative DM. The proof of the following theorem is included in the Appendix \ref{APP:the5}.

\begin{theorem}\label{the:Approx-MK}
Consider the MDP for the representative DM under Assumption \ref{assump:1}. Suppose further that 
\begin{itemize}
    \item [(i)] $\mu_{0}$ is non-atomic,
    \item [(ii)] $\mathcal{T}(\cdot|x_0^R,u_{0}^{R},\mu_{0})$ is non-atomic for every $x_0^R,u_{0}^{R}$ and $\mu_{0}$.
\end{itemize}
Let an optimal solution for ${\sf MDP_n}$ attain the optimal performance $\widehat{J}_{0,T}^{n}$ for the finite horizon problem and $\widehat{J}_{\infty}^{n}$ for the infinite horizon problem. Then, 
\begin{align}\label{the:Ap-main}
    &\lim_{n\to \infty}|\widehat{J}_{0,T}^{n}(\mu_0)-J_{0,T}(\mu_0)|=0,\\
    &\lim_{n\to \infty}|\widehat{J}_{\infty}^{n}(\mu_0)-J_{\infty}(\mu_0)|=0\nonumber,
\end{align}
where $J_{0,T}(\mu_0)$ and $J_{\infty}(\mu_0)$ are the optimal performance for the finite horizon and infinite horizon MDP for the representative DM, respectively.
\end{theorem}

In the above, we assume that the cost and transition kernels are known to DMs. In learning scenarios, DMs do not have access to this knowledge and follow a learning algorithm balancing exploration and exploitation that can be implemented either independently (see e.g., \cite{yongacoglu2022independent}) or in a centralized fashion (see e.g., \cite{anahtarci2023learning, carmona2023model}). We now provide a brief discussion on learning for MF teams.

In \cite{yongacoglu2022independent}, a win-stay, lose-shift algorithm has been introduced that guarantees asymptotic convergence to an approximate subjective equilibrium under MF sharing information structure for MF games with a finite number of DMs. In addition, if the obtained policy is symmetric, then it constitutes an approximate Nash equilibrium. Utilizing this result in our setting, we can show that the algorithm  \cite{yongacoglu2022independent} converges to an approximate subjective person-by-person optimal solution. Also, it converges to an approximate person-by-person optimal solution if it is symmetric. This, however, does not imply global optimality since the limiting policy will possibly be subjective and/or locally optimal. Devising an independent learner algorithm that converges to an approximate globally optimal solution remains open and left for future work. On the other hand, thanks to Theorem \ref{the:3}, the representative DM model can be utilized to approximate a symmetric globally optimal solution by considering offline learning algorithms, assuming that DMs have access to an oracle or a simulator. In this direction, \cite{anahtarci2023learning, carmona2023model,kara2023q,subramanian2023mean} might apply to this setting that can be implemented in a centralized fashion by forcing DMs to apply symmetric policies.

\section{Conclusion.}
We have studied decentralized stochastic exchangeable teams with MF sharing information comprising a finite number of DMs, as well as their MF limits involving infinite numbers of DMs. For finite population exchangeable teams, we have established the existence of a randomized optimal policy that is exchangeable and Markovian via value iterations for an equivalent measure-valued controlled MDP. We then established the existence of a symmetric, independent, decentralized optimal randomized policy for the infinite population problem and proved the optimality of the limiting measure-valued MDP for the representative DM. Also, we have shown the approximate optimality of symmetric, independent, decentralized optimal randomized policies for the corresponding finite-population team with many DMs under the centralized information structure. Finally, we have provided a verification theorem and a finite MDP to approximate a solution of the MDP for the representative DM.

\appendix

\section{Proof of Section \ref{sec:Randomization}.}

\subsection{Proof of Lemma \ref{lem:sym}.}\label{APP:Lem-sym}

In the following, we only prove \eqref{eq:lem-eq2}. The proof of \eqref{eq:lem-eq1} follows from an analogous argument.
Fix $\mu_{t}^{N}$. Suppose that $x_{t}^{1:N}$ and $\widehat{x}_{t}^{1:N}$ have the same empirical measure  $\mu_{t}^{N}$. Then, we must have that $\widehat{x}_{t}^{i}=x_{t}^{\sigma(i)}$ almost surely for $i=1,\ldots, N$ and for a permutation $\sigma\in S_{N}$. Denote the set of all $x_{t}^{1:N}$ with empirical measure $\mu_{t}^{N}$ by $A(\mu_{t}^{N})$. Similarly, fix $\theta_{t}^{N}$. Suppose that $(x_{t}^{1:N},u_{t}^{1:N})$ and $(\widehat{x}_{t}^{1:N},\widehat{u}_{t}^{1:N})$ have the same joint state-action empirical measure  $\theta_{t}^{N}$. Then, we must have that $\widehat{x}_{t}^{i}=x_{t}^{\sigma(i)}$ and $\widehat{u}_{t}^{i}=u_{t}^{\sigma(i)}$ almost surely for $i=1,\ldots, N$ and for a permutation $\sigma\in S_{N}$. Denote the set of all $u_{t}^{1:N}$ with  $\theta_{t}^{N}$ by $A(\theta_{t}^{N})$.

In the following, we show that the following chain of equalities holds:
\begin{align}
    Pr\left((\mu_{t+1}^{N},x_{t+1}^{i})\in \cdot | x_{t}^{i},\overline{\pi}_{t}^N,\mu_{t}^{N}\right)
&=Pr\left((\mu_{t+1}^{N},x_{t+1}^{i})\in \cdot | x_{t}^{i},\overline{\pi}_{t}^N,\mu^{N}_t, \theta^N_{t}\right)\label{eq:e2}\\
    &=Pr\left((\mu_{t+1}^{N},x_{t+1}^{i})\in \cdot | x_{t}^{i},\overline{\pi}_{t}^N,\mu^{N}_t, \theta^{N}_t, x_{t}^{1:N},u_{t}^{1:N}\right)\label{eq:e3}\\
    &=Pr\left((\mu_{t+1}^{N},x_{t+1}^{i})\in \cdot | x_{t}^{i},\overline{\pi}_{t}^N,\mu^{N}_{0:t}, \theta^{N}_{0:t}, x_{0:t}^{1:N},u_{0:t}^{1:N}\right)\label{eq:e5}
\end{align}
for any symmetric independent policy $\overline{\pi}_{t}^N$ for all DMs.
Given $\mu_{t}^{N}$, any  $x_{t}^{1:N}\in A(\mu_{t}^{N})$ are permutation of some states $x_{t}^{1:N}$. Since policies are symmetric for any $x_{t}^{1:N}\in A(\mu_{t}^{N})$, we have 
\begin{align*}
   \frac{1}{N}\sum_{i=1}^{N}\delta_{(x_{t}^{\sigma(i)},\widehat{u}_{t}^{i})}(\cdot)
    &= \frac{1}{N}\sum_{i=1}^{N}\delta_{(x_{t}^{i},{u}_{t}^{i})}(\cdot):=\theta_{t}^N
\end{align*}
for any permutation $\sigma\in S_N$, where ${u}_{t}^{i}\sim \overline{\pi}_{t}(\cdot|h_{t}^i)$ and $\widehat{u}_{t}^{i}\sim \overline{\pi}_{t}(\cdot|h_{t}^{\sigma(i)})$ with $h_{t}^i:=(x_{0:t}^{i},\mu_{0:t}^N)$ and $h_{t}^{\sigma(i)}:=(x_{0:t}^{\sigma(i)},\mu_{0:t}^N)$ for all $t\geq 0$.
This implies that for all $x_{k}^{1:N}\in A(\mu_{k}^{N})$ with $k=0,\ldots, t$, all possible $u^{1:N}_{t}$ are such that the joint state-action empirical measure is $\theta_{t}^{N}$ with marginals on states  $\mu_{t}^{N}$. This implies that $u_{t}^{1:N}\in A(\theta_{t}^{N})$, and hence, \eqref{eq:e2} holds. All $x_{t}^{1:N}\in A(\mu_{t}^{N})$ and $u_{t}^{1:N}\in A(\theta_{t}^{N})$ are permutation of some state and action pairs $(x_{t}^{1:N},u_{t}^{1:N})$, and hence, using the fact that dynamics is symmetric (i.e., $f$ is the for  each DM in \eqref{eq:dyn}) and $w_{t}^{1:N}$ in \eqref{eq:dyn} are i.i.d., we have for any $\sigma\in S_N$
    \begin{align*}
        Pr(x_{t+1}^{i} \in \cdot|x_{t}^{\sigma(1):\sigma(N)},u_{t}^{\sigma(1):\sigma(N)}, \mu_{t}^{N}, \overline{\pi}_{t}^N)
        &=Pr(x_{t+1}^{\sigma(i)} \in \cdot|x_{t}^{1:N},u_{t}^{1:N}, \mu_{t}^{N}, \overline{\pi}_{t}^N).
    \end{align*}
This implies that 
\begin{align*}
        Pr(x_{t+1}^{1:N} \in \cdot|x_{t}^{\sigma(1):\sigma(N)},u_{t}^{\sigma(1):\sigma(N)}, \mu_{t}^{N}, \overline{\pi}_{t}^N)
        &=\prod_{i=1}^{N}Pr(x_{t+1}^{i} \in \cdot|x_{t}^{\sigma(1):\sigma(N)},u_{t}^{\sigma(1):\sigma(N)}, \mu_{t}^{N}, \overline{\pi}_{t}^N)\\
        &=\prod_{i=1}^{N}Pr(x_{t+1}^{\sigma(i)} \in \cdot|x_{t}^{1:N},u_{t}^{1:N}, \mu_{t}^{N}, \overline{\pi}_{t}^N)\\
&=Pr(x_{t+1}^{\sigma(1):\sigma(N)} \in \cdot|x_{t}^{1:N},u_{t}^{1:N}, \mu_{t}^{N}, \overline{\pi}_{t}^N).
\end{align*}
Hence, the set of all $x_{t+1}^{1:N}$ are conditionally permutations of each other in distribution for all $x_{t}^{1:N}\in A(\mu_{t}^{N})$ and $u_{t}^{1:N}\in A(\theta_{t}^{N})$. In other words, if $x^{1:N}_{t+1}$ and $\widehat{x}^{1:N}_{t+1}$ are two possible states, then $\mathcal{L}(\widehat{x}^{1:N}_{t+1})=\mathcal{L}({x}^{\sigma(1):\sigma(N)}_{t+1})$ for a permutation $\sigma\in S_N$. This implies that $Pr(\mu_{t+1}^{N}\in \cdot | x_{t}^{1:N}, u_{t}^{1:N},\mu^{N}_t, \theta^N_{t})$ remains the same for all $x_{t}^{1:N}\in A(\mu_{t}^{N})$ and $u_{t}^{1:N}\in A(\theta_{t}^{N})$. Since for each DM$^{i}$, $x_{t}^{i}$ is independent of $x_{t}^{-i},u_{t}^{-i}$ given $\mu_{t}^{N},x_{t}^{i},\overline{\pi}_{t}^N$, we get \eqref{eq:e3}. Equality \eqref{eq:e5} follows from the fact that $x_{t+1}^{1:N}$ is independent of $\overline{\pi}_{t}^N,\mu^{N}_{0:t}, \theta^{N}_{0:t}, x_{0:t}^{1:N},u_{0:t}^{1:N}$ given $x_{t}^{1:N},u_{t}^{1:N}$. This completes the proof.

\section{Proofs of Section \ref{sec:Existence}.}

\subsection{Proof of Theorem \ref{the:2}.}\label{APP-t2}
Part (i): The proof proceeds in two steps:
\begin{itemize}[wide]
    \item {\bf Step 1: Exchangeability property of an optimal solution.} We start at $t=T-1$, suppose that $g_{T-1}^{\star}$ is optimal, i.e.,
\begin{align*}
   J^{N}_{T-1,T}(\mu^N):=\inf_{\theta^N\in U(\mu^N)}\tilde{c}(\mu, \theta^N)= \tilde{c}(\mu^N, g_{T-1}^{\star}(\mu^N)).
\end{align*}
This induces $g_{T-1}^{\star}(\mu^N)=\theta^{N\star}_{T-1}(x^{1:N}_{T-1},u^{1:N}_{T-1})$. In particular, $\theta^{N\star}_{T-1}$ induces a joint distribution on $x^{1:N}_{T-1},u^{1:N}_{T-1}$ with marginals on $x^{1:N}_{T-1}$ coincides with $\mu^N$. Denote this joint distribution by $P^{\star}_{T-1}$. We note that $P^{\star}_{T-1}$ is induced by a Markovian randomized policy $\pi_{T-1}^{\star}$, i.e.,
\begin{align*}
    P^{\star}_{T-1}(\cdot)=\int_{\cdot}\pi_{T-1}^{\star}(du_{T-1}^{1:N}|x_{T-1}^{1:N})Pr(dx_{T-1}^{1:N}).
\end{align*}

Denote the permutation of $P^{\star}_{T-1}$ by $P^{\sigma\star}_{T-1}$. This permutation does not change $\mu^N$ and $\theta^{N\star}_{T-1}$; however, $P^{\star}_{T-1}$ does not need to be exchangeable. Construct a new joint distribution on  $x^{1:N}_{T-1},u^{1:N}_{T-1}$ by averaging over all permutation of $P^{\star}_{T-1}$ as follows
\begin{align*}
\widehat{P}^{\star}_{T-1}=\frac{1}{|S_{N}|}\sum_{\sigma\in S_N} P^{\sigma\star}_{T-1}.
\end{align*}
We have that $\widehat{P}^{\star}_{T-1}$ is exchangeable.  We have
\begin{align}
    &\int \widetilde{c}(\mu(x^{1:N}_{T-1}), \theta^N(x^{1:N}_{T-1},u^{1:N}_{T-1})) \widehat{P}^{\star}_{T-1}(dx^{1:N}_{T-1},du^{1:N}_{T-1})\label{eq:1-1}\\
    &=\frac{1}{|S_{N}|}\sum_{\sigma\in S_N}\int \tilde{c}(\mu(x^{1:N}_{T-1}), \theta^N(x^{1:N}_{T-1},u^{1:N}_{T-1}))  {P}^{\sigma\star}_{T-1}(dx^{1:N}_{T-1},du^{1:N}_{T-1})\label{eq:lin1}\\
    &=\frac{1}{|S_{N}|}\sum_{\sigma\in S_N}\int \tilde{c}(\mu(x^{\sigma(1):\sigma(N)}_{T-1}), \theta^N(x^{\sigma(1):\sigma(N)},u^{\sigma(1):\sigma(N)}))  {P}^{\star}_{T-1}(dx^{1:N}_{T-1},du^{1:N}_{T-1})\\
    &=\frac{1}{|S_{N}|}\sum_{\sigma\in S_N}\int \tilde{c}(\mu(x^{1:N}_{T-1}), \theta^N(x^{1:N}_{T-1},u^{1:N}_{T-1})) {P}^{\star}_{T-1}(dx^{1:N}_{T-1},du^{1:N}_{T-1})\label{eq:emp-does-not-change}\\
    &=\int \tilde{c}(\mu(x^{1:N}_{T-1}), \theta^N(x^{1:N}_{T-1},u^{1:N}_{T-1})) {P}^{\star}_{T-1}(dx^{1:N}_{T-1},du^{1:N}_{T-1})\\
    &=\tilde{c}(\mu,g_{T-1}^{\star}(\mu^N))\label{eq:def-P^*},
\end{align}
where \eqref{eq:lin1} follows from the definition of $\widehat{P}^{\star}_{T-1}$ and linearity of expectation in $\widehat{P}^{\star}_{T-1}$. Equality \eqref{eq:emp-does-not-change} follows from the fact that the permutation of $P^{\star}_{T-1}$, denoted by $P^{\sigma\star}_{T-1}$ does not change empirical measures of states and actions, and \eqref{eq:def-P^*} follows from the definition of $P^{\star}_{T-1}$ since it is induced by $\theta^{N\star}_{T-1}$, which is determined by $g^{\star}_{T-1}$. This implies that without any loss, we can assume that at $t=T-1$, the optimal joint distribution of $x^{1:N}_{T-1},u^{1:N}_{T-1}$ is exchangeable, where the actions are induced by a permutation invariant Markov randomized policy.

Now, we consider $t=T-2$. Suppose that $g_{T-2}^{\star}$ is optimal, i.e.,
\begin{align}
    J^{N}_{T-2,T}(\mu^N)&:=\inf_{\theta^N\in U(\mu^N)}\left\{\tilde{c}(\mu, \theta^N)+\beta \int J^{N}_{T-1,T}(\mu^N_{+}) \eta^N(d\mu^N_{+}|\mu^N, \theta^N)\right\}\nonumber\\
    &=\tilde{c}(\mu, g_{T-2}^{\star}(\mu^N))+\beta \int J^{N}_{T-1,T}(\mu^N_{+}) \eta^N(d\mu^N_{+}|\mu^N, g_{T-2}^{\star}(\mu^N))\label{eq:DP-T-2}
\end{align}
This induces $g_{T-2}^{\star}(\mu^N)=\theta^{N\star}_{T-2}(x^{1:N}_{T-2},u^{1:N}_{T-2})$. In particular, $\theta^{N\star}_{T-2}$ induces a joint distribution on $x^{1:N}_{T-2},u^{1:N}_{T-2}$ with marginals on $x^{1:N}_{T-2}$ coincides with $\mu^N$. Denote this joint distribution by $P^{\star}_{T-2}$. Note the permutation of $P^{\star}_{T-2}$, denoted by $P^{\sigma\star}_{T-2}$ does not change $\mu^N$ and $\theta^{N\star}_{T-2}$; however, $P^{\star}_{T-2}$ does not need to be exchangeable. Construct a new joint distribution on  $x^{1:N}_{T-2},u^{1:N}_{T-2}$ by averaging over all permutation of $P^{\star}_{T-2}$ as follows
\begin{align*}
\widehat{P}^{\star}_{T-1}=\frac{1}{|S_{N}|}\sum_{\sigma\in S_N} P^{\sigma\star}_{T-2}.
\end{align*}
We have that $\widehat{P}^{\star}_{T-2}$ is exchangeable. An analogous argument as used in \eqref{eq:1-1}---\eqref{eq:def-P^*} yields that
\begin{align}\label{eq:47}
    &\int \widetilde{c}(\mu^N(x^{1:N}_{T-2}), \theta^N(x^{1:N}_{T-2},u^{1:N}_{T-2})) \widehat{P}^{\star}_{T-2}(dx^{1:N}_{T-2},du^{1:N}_{T-2})=\widetilde{c}(\mu^N,g_{T-2}^{\star}(\mu^N)).
\end{align}
Similarly, we can show that
\begin{align*}
    &\int \int J^{N}_{T-1,T}(\mu^N_{+}) \eta^N\left(d\mu^N_{+}|\mu^N(x_{T-2}^{1:N}), \theta^N(x_{T-2}^{1:N},u_{T-2}^{1:N})\right)\widehat{P}^{\star}_{T-2}(dx^{1:N}_{T-2},du^{1:N}_{T-2})\\
    &=\frac{1}{|S_{N}|}\sum_{\sigma\in S_N} \int \int J^{N}_{T-1,T}(\mu^N_{+}) \eta^N\left(d\mu^N_{+}|\mu^N(x_{T-2}^{1:N}), \theta^N(x_{T-2}^{1:N},u_{T-2}^{1:N})\right){P}^{\sigma\star}_{T-2}(dx^{1:N}_{T-2},du^{1:N}_{T-2})\\
    &=\int J^{N}_{T-1,T}(\mu^N_{+}) \eta^N(d\mu^N_{+}|\mu^N, g_{T-2}^{\star}(\mu^N)).
\end{align*}
This together with \eqref{eq:DP-T-2} and \eqref{eq:47}  implies that we can assume that at $t=T-2$, the optimal joint distribution of $x^{1:N}_{T-2},u^{1:N}_{T-2}$ is exchangeable without any loss of optimality, where the actions are induced by a permutation invariant Markov randomized policy. Sequentially, we can show that for any $t=0,\ldots, T-1$, an optimal joint distribution of $x^{1:N}_{T-2},u^{1:N}_{T-2}$ is exchangeable, where the actions are induced by a permutation invariant Markov randomized policy. This implies that if the value iteration \eqref{eq:value-iter-T} admits a solution $\pmb{g_{T}^{\star}}$, then without any loss, we can restrict the search for an optimal strategic measure for \PNT\ to exchangeable ones $\LEXN$ induced by exchangeable Markov randomized policies in $\PEXNM$.

\item {\bf Step 2: Weak Feller property of MDP and continuity of the cost.} Using Assumption \ref{assump:1}, by the classical measurable selection theorem in the stochastic control literature, following the measure-valued MDP formulation, we can show that there exists $\pmb{g_{T}^{\star}}\in \prod_{t=0}^{T-1}\mathbb{G}_t$ for \PNT.

Since \(N\) is fixed, if $\theta_k^N$ converges weakly to \(\theta^N\) as $k\to \infty$, then (after possibly relabeling) we get
$(x_k^i,u_k^i)$ converges to $(x^i,u^i)$ as $k\to \infty$,
for  $i=1,\dots,N$. For each \(k\), let \(x^1_{+,k},\dots,x^N_{+,k}\) are independent with
$x^i_{+,k}\sim \mathcal T(\cdot\mid x_k^i,u_k^i,\mu_k^N)$,
and \((x^1_{+},\dots,x^N_{+})\) are independent with
$x^i_{+}\sim \mathcal T(\cdot\mid x^i,u^i,\mu^N)$.
Define $\mu_{+,k}^{N}$ and $\mu_{+}^{N}$ as the empirical measures $x^{1:N}_{+,k}$ and $x^{1:N}_{+}$, respectively. We show that
$\eta^N(\cdot\mid \mu_k^N,\theta_k^N)=\mathcal L(\mu_{+,k}^{N})$ converges to $\eta^N(\cdot\mid \mu^N,\theta^N)=\mathcal L(\mu_{+}^{N})$. 

Let \(\Phi\in C_b(\mathcal P(\mathbb X))\), and define
$G_N(y^{1:N}):=
\Phi(\frac1N\sum_{i=1}^N \delta_{y^i})$
for $y^i=x^i_{+,k}$ or $y^i=x^i_{+}$.
Then, \(G_N\in C_b(\mathbb X^N)\) since \(\Phi\in C_b(\mathcal P(\mathbb X))\) and the empirical measure is continuous in $y^{1:N}$. By Assumption \ref{assump:1}(iii), we get 
\begin{align*}
   \lim_{k\to \infty} \int \Phi(\nu)\,\eta^N(d\nu\mid \mu_k^N,\theta_k^N)&=
\lim_{k\to \infty}\int G_N(x^{1:N}_{+})\prod_{i=1}^N \mathcal T(dx^i_{+}\mid x_k^i,u_k^i,\mu_k^N)\\
&=\int G_N(x^{1:N}_{+})\prod_{i=1}^N \mathcal T(dx^i_{+}\mid x^i,u^i,\mu^N)\\
&=\int \Phi(\nu)\,\eta^N(d\nu\mid \mu^N,\theta^N).
\end{align*}

Next, we show that the cost is also continuous in $\theta_t^N$ and $\mu_t^N$. Suppose that $\theta_{t,k}^N$ and $\mu_{t,k}^N$ converge weakly to $\theta_t^N$ and $\mu_t^N$ as $k\to \infty$. Under Assumption \ref{assump:1}(ii), we have
\begin{align*}
    \lim_{k\to \infty} \tilde{c}(\mu_{t,k}^N,\theta_{t,k}^N)&=\lim_{k\to \infty} \int c(x_t,u_t,\mu_{t,k}^N)\theta_{t,k}^N(dx_t,du_t)\\
    &=\int c(x_t,u_t,\mu_{t}^N)\theta_{t}^N(dx_t,du_t).
\end{align*}
This completes the proof of part (i) since the measure-valued MDP has a weak Feller property with a non-negative continuous cost, and hence, $\pmb{g^{\star}_T}$ exists.

 Part (ii): Using the monotone convergence theorem and the analysis for the classical MDP problem, the limit of $J^N_{\infty}(\mu):=\lim\limits_{T\to \infty}J^{N}_{t,T}(\mu)$ exists using a fixed-point argument. Hence, there exists $\pmb{g^{\star}}\in \prod_{t=0}^{\infty}\mathbb{G}_t$ for \PN\ that is Markovian and stationary, satisfying the value iterations \eqref{eq:value-iter-inf}. Using an analogous argument as that used in the proof of part (i), we can conclude that the optimal solution for \PN\ is induced by an exchangeable strategic measure $\LEXN$ with a stationary Markov randomized policy in $\PEXNS$. This completes the proof of part (ii).
\end{itemize}

\subsection{Proof of Theorem \ref{the:3}.}\label{APP-t3}

To prove Theorem \ref{the:3}, we use the following two key results by Aldous \cite[Proposition 7.20]{aldous2006ecole} and  Diaconis and Friedman \cite[Theorem 13]{diaconis1980finite}, which we recall in the following.

\begin{theorem}\cite[Proposition 7.20]{aldous2006ecole}\label{the:ald}
Let $X:=(X_{1},X_{2},\dots)$ be an infinitely exchangeable sequence of random variables taking values in a Polish space $\mathbb{X}$. Suppose that $X$ is directed by a random measure $\alpha:\Omega \to \mathcal{P}(\X)$, i.e.,
\begin{align}\label{mixture-iid}
    P(X_{1}\in A^1, X_{2}\in A^2, \ldots)= \int \prod_{i=1}^{\infty}\alpha(A^i) \eta(d\alpha) 
\end{align}
In other words, conditioned on $\alpha$, $(X_{1},X_{2},\dots)$ are i.i.d. random variables. 
Suppose that 
\begin{enumerate}
\item either for each $n$, $X^{(n)}=(X_{1}^{(n)},X_{2}^{(n)},\dots)$  is infinitely exchangeable directed by $\alpha_{n}$,
\item or $X^{(n)}=(X_{1}^{(n)},\dots, X_{n}^{(n)})$  is $n$-exchangeable with empirical measure $\alpha_{n}$.
\end{enumerate}
Then, $X^{(n)}$ converges in distribution\footnote{By convergence in distribution to an infinitely exchangeable sequence, we mean the following: $X^{(n)} \xrightarrow[n \to \infty]{{\sf d}} X$ if and only if $(X_{1}^{(n)},\dots, X_{m}^{(n)})\xrightarrow[n \to \infty]{{\sf d}}(X_{1},\dots, X_{m})$ for each $m\geq 1$ \cite[page 55]{aldous2006ecole}.} to $X$ $\big(X^{(n)}\xrightarrow[n \to \infty]{{\sf d}}X\big)$  if and only if $\alpha_{n}\xrightarrow[n \to \infty]{{\sf d}}\alpha$.
\end{theorem}

\begin{theorem}\cite[Theorem 13]{diaconis1980finite}\label{the:dia}
Let $Y=(Y_{1},\dots,Y_{n})$ be an $n$-exchangeable and $Z=(Z_{1},Z_{2},\dots)$ be an infinitely exchangeable sequence of random variables with $\mathcal{L}(Z_{1},\dots,Z_{k})=\mathcal{L}(Y_{I_1},\dots,Y_{I_k})$ for all $k\geq 1$, where the indices $(I_{1}, I_{2}, \dots)$ are i.i.d. random variables with the uniform distribution on the set $\{1,\dots,n\}$. Then,  for all $m=1,\dots, n$,
\begin{flalign*}
\left|\left|\mathcal{L}(Y_{1},\dots,Y_{m})-\mathcal{L}(Z_{1},\dots,Z_{m})\right|\right|_{\sf TV}\leq \frac{m(m-1)}{2n},
\end{flalign*}
where $\mathcal{L}(\cdot)$ denotes the law of random variables and $||\cdot||_{\sf TV}$ is the total variation norm.
\end{theorem}

\begin{proof}[Proof of Theorem \ref{the:3}.]
We proceed in five steps:
\begin{itemize}[wide]
    \item {\bf Step 1: Convergence of optimal strategic measures.}
    {Suppose that $\theta^{N\star}_t$ with marginals $\mu^{N\star}_t$ (which are induced by an exchangeable strategic measure ${P_{\pi}^{N\star}}$) solve the value iteration equations \eqref{eq:value-iter-T} for every $t\geq 0$. Denote the states and actions under ${P_{\pi}^{N\star}}$ by $({x^{1\star:N\star}_{t,N}},{u^{1\star:N\star}_{t,N}})$ for $t\geq 0$ that is $N$-exchangeable. For each $t\geq 0$, by Theorem \ref{the:dia}, there exists an infinitely exchangeable sequence $({\widehat{x}^{1\star:N\star}_{t,N}},{\widehat{u}^{1\star:N\star}_{t,N}})$ such that
\begin{flalign}\label{eq:frid-policy}
\left\|\mathcal{L}({x^{1\star:m\star}_{t,N}},{u^{1\star:m\star}_{t,N}})-\mathcal{L}({\widehat{x}^{1\star:m\star}_{t,N}},{\widehat{u}^{1\star:m\star}_{t,N}})\right\|_{\sf TV}\leq \frac{m(m-1)}{2N}
\end{flalign}  
for all $m=1,\ldots, N$. 

Since $\X$ and $\U$ are compact, the joint probability measures on $\X$ and $\U$ are tight and hence (weakly) relatively compact. Thus, for every $t\geq 0$, $(\underline{{\widehat{x}_{t,N}^{\star}}},\underline{{\widehat{u}_{t,N}^{\star}}}):=({\widehat{x}^{1\star:\infty\star}_{t,N}},{\widehat{u}^{1\star:\infty\star}_{t,N}})$ admits a subsequence  $(\underline{{\widehat{x}_{t,n}^{\star}}},\underline{{\widehat{u}_{t,n}^{\star}}})$ that converges in distribution to $(\underline{{{x}_t^{\star}}},\underline{{{u}_t^{\star}}})$, that is 
 \begin{align*}
\left({\widehat{x}_{t,n}^{1\star:m\star}},{\widehat{u}_{t,n}^{1\star:m\star}}\right)    \xrightarrow[n \to \infty]{{\sf d}} \left({{x}^{1\star:m\star}_{t}},{{u}^{1\star:m\star}_{t}}\right) \quad \forall m\geq 1.
 \end{align*}
 The limiting process $\{(\underline{{{x}^{\star}_{t}}},\underline{{{u}^{\star}_{t}}})\}_{t\geq0}$ is also infinitely-exchangeable. Since for any $t\geq 0$,  $({\widehat{x}^{1\star:m\star}_{t,n}},{\widehat{u}^{1\star:m\star}_{t,n}})$ converges weakly, using \eqref{eq:frid-policy}, any finite marginal of ${{P}_{\pi}^{n\star}}$ on states and actions $({{x}^{1\star:m\star}_{t,n}},{{u}^{1\star:m\star}_{t,n}})$ also converges weakly to the marginals ${{P}_{\pi}^{\star}}$ on states and actions $({{x}^{1\star:m\star}_{t}},{{u}^{1\star:m\star}_{t}})$ as $n\to \infty$ for every $m\geq 1$. This implies that $n$-exchangeable random variables $({{x}^{1\star:n\star}_{t,n}},{{u}^{1\star:n\star}_{t,n}})$ converges in distribution to the infinitely exchangeable random variables $(\underline{{{x}_t^{\star}}},\underline{{{u}_t^{\star}}})=({{x}^{1\star:\infty\star}_t},{{u}^{1\star:\infty\star}_t})$ for each $t\geq 0$.

Using Theorem \ref{the:ald}, this implies that the sequence of empirical measures $\{{\theta^{n\star}_t}\}_n$ of $\{({{x}^{1\star:n\star}_{t,n}},{{u}^{1\star:n\star}_{t,n}})\}_n$ together with their marginals $\{{\mu^{n\star}_t}\}_n$ converges in distribution to the directing measure ${\theta^{\star}_t}$ with its marginal ${\mu^{\star}_t}$ of $(\underline{{{x}_t^{\star}}},\underline{{{u}_t^{\star}}})$ as $n\to \infty$ for each $t\geq 0$. This implies that $\{{P_{\pi}^{N\star}}\}_n$ converges to ${{P}_{\pi}^{\star}}$ weakly (for any finite marginals) as $n\to \infty$. Additionally, it yields that the directing measure of the limiting state and action process is the weak limit of the empirical measure of the optimal state and action process of the $N$-agent problem attained via the value iteration. Consequently, this implies that the weak limit leads to a control policy that is computable by only knowing the realization of the state and the MF. The limiting strategic measure is conditionally i.i.d. given the weak limit $\pmb{\theta^{\star}}$ (with marginal $\pmb{\mu^{\star}}$) of the empirical measure on the induced optimal states and actions for \PNT\ and \PN.

As a result, for each time, using \eqref{mixture-iid}, we can write  at each given time $t\geq 0$ for the limiting strategic measure, the joint distribution on states and actions conditioned on the directing measure $\theta_t^{\star}$ with marginal $\mu_t^{\star}$ as 
\begin{align}\label{eq:rep-policy}
Pr((\underline{u}_t,\underline{x}_t)\in \cdot |\theta_t^{\star})
 &=\prod_{i=1}^{\infty}Pr((u^i_t,x_t^i)\in \cdot|\theta_{t}^{\star})\\
 &=\prod_{i=1}^{\infty}\theta_{t}^{\star}((u^i_t,x_t^i)\in \cdot)\nonumber\\
 &=\prod_{i=1}^{\infty}\pi_t^{\star}(u^i_t \in \cdot|x_t^i) \mu_t^{\star}(\cdot)\nonumber
\end{align}
for some stochastic kernel $\pi_t^{\star}$ which will be selected using $\mu_t^{\star}$ (under MF information sharing). 

Next, we show that under any infinite exchangeable policy, the limiting problem constitutes an MDP. Since $(x_t^{1:\infty},u_t^{1:\infty})$ are infinitely exchangeable with the directing measure $\theta_t$, conditional on $\theta_t$, we have $(x_t^{1:\infty},u_t^{1:\infty})$ are i.i.d. Also, conditioned on $\mu_t$, $(x_t^{1:\infty})$ are i.i.d. Since 
$x_{t+1}^{i}\sim \mathcal{T}(\cdot|x_t^i,u_t^i,\mu_t)$, $x_{t+1}^{1:\infty}$ are i.i.d. conditioned on realization of $\theta_t,\mu_t$ with 
\begin{align*}
F(\mu_t,\theta_t)(\mathbb{A}):=P\!\left(x^{i}_{t+1}\in \mathbb{A} \mid \mu_t,\theta_t\right)
&=
\int_{\mathbb{A}}
\mathcal{T}\!\left(dx^{i}_{t+1}\mid x_t^i,u_t^i,\mu_t\right)\,
\theta_t(dx_t^i,du_t^i).
\end{align*}
for every Borel set $\mathbb{A}$ subset of $\X$.
Thus, $\mathcal{L}(x_{t+1}^{i}|\theta_t)=F(\mu_t,\theta_t)$.

Let $\mu_{t+1\mid N}$ and $\theta_{t+1\mid N}$ be the empirical measures of infinitely exchangeable sequence of $(x_{t+1}^{1:\infty})$ and $(x_{t+1}^{1:\infty},u_{t+1}^{1:\infty})$, respectively. Since conditioned on $\theta_t$, $x_{t+1}^{1:\infty}$ are i.i.d. with distribution $F(\mu_t,\theta_t)$, by strong law of large numbers, we get $\mu_{t+1\mid N}$ converges weakly to $F(\mu_t,\theta_t)$. Hence, as $N\to \infty$, we can define the transition kernel of the MDP as
\begin{align*}
\eta(\,\mathbb{B} \mid \mu_t,\theta_t):=P\!\left(\mu_{t+1} \in \mathbb{B} \mid \mu_t,\theta_t\right)
&=
\delta_{\{F(\mu_t,\theta_t)\in \mathbb{B}\}}
\end{align*}
for any Borel set $\mathbb{B}$ subset of $\mathcal{P}(\X)$. Next, we obtain the running cost for this MDP. We have 
\begin{align*}
\limsup_{N\to\infty}\frac{1}{N}\sum_{i=1}^N
c\!\left(x_t^i,u_t^i,\mu_{t\mid N}\right)
&=
\limsup_{N\to\infty}
\int c\!\left(x_t,u_t,\mu_{t\mid N}\right)\,
\theta_{t\mid N}(dx,du),
\end{align*}
where $\mu_{t\mid N}$ and $\theta_{t\mid N}$ is the empirical measure of infinitely exchangeable sequence of $(x_t^{1:\infty})$ and $(x_t^{1:\infty},u_t^{1:\infty})$. Since $(x_t^{1:\infty},u_t^{1:\infty})$ are infinitely exchangeable, conditioned on $\theta_t$, they are independent and identical, and hence, by the law of large numbers, $\theta_{t\mid N}$ and $\mu_{t\mid N}$ converges weakly to $\theta_t$ and $\mu_t$, respectively. This implies that the running cost of the MDP is given by 
\begin{align*}
\tilde c(\mu_t,\theta_t):=\int c\!\left(x_t,u_t,\mu_t\right)\,
\theta_t(dx,du)=\lim_{N\to \infty}\int c\!\left(x_t,u_t,\mu_{t\mid N}\right)\,
\theta_{t\mid N}(dx,du).
\end{align*}
This implies the measure-valued MDP formulation under infinite exchangeable policies with the well-defined stochastic kernel $\eta$ and the running cost $\tilde c$.
\item {\bf Step 2: Convergence of the $N$-agent value function.} Let $J^{N}_{0,T}$ be the value function under an optimal strategic measure $\{{P_{\pi}^{N\star}}\}_n$ for \PNT\ and $J_{0,T}$ be the value of the limiting strategic measure ${{P}_{\pi}^{\star}}$. We have
\begin{align*}
    \inf_{\underline{\pmb{{\gamma}}}\in \prod_{i=1}^{\infty}\pmb{\Gamma^{i}}}
J(\underline{\pmb{{\gamma}}})&\geq \limsup_{N\to \infty}\inf_{{P_{\pi}^N}\in \LN} J_{T}^{N}({P_{\pi}^N})\\
    &\geq \lim_{n\to \infty} J_{0,T}^{n}\\
    &= J_{0,T},
\end{align*}
where the first inequality follows from exchanging limsup and inf and the fact that deterministic policies are a strict subset of randomized policies, and the second inequality follows from the fact that the limsup is the greatest subsequential limit. In the following, we justify the last equality. To this end, we first show that the stochastic kernel $\eta^N(\cdot\mid \mu^{N}_t,\theta^{N}_t)$ of finite agent measure-valued MDP  converges continuously (in \(\mathcal P(\mathbb X)\)) to the stochastic kernel $\eta(\cdot\mid \mu_t,\theta_t)$ of the limiting measure-valued MDP as $N\to \infty$ for all $t\geq 0$. We equip \(\mathcal P(\mathbb X)\) with the following metric \(\rho\) generating the weak topology
\begin{align*}
\rho(\mu,\nu)=\sum_{m=1}^{\infty}2^{-(m+1)}\bigl|\mu(f_m)-\nu(f_m)\bigr|,
\end{align*}
where \(\{f_m\}\subset C_b(\mathbb X)\) is chosen so that \(\|f_m\|_\infty\le 1\) for all \(m\). We prove that if $(\mu_t^N,\theta_t^N)$ converges weakly to $(\mu_t,\theta_t)$ as $N\to \infty$, then $\mathcal{L}(\mu_{t+1}^N \mid \mu_t^N,\theta_t^N)(\cdot)$ converges to 
$\delta_{F(\mu_t,\theta_t)}(\cdot)$ as $N\to \infty$.

For any \(f\in C_b(\mathbb X)\), define
$\mathcal T f(x,u,\mu)
:=\int_{\mathbb X} f(y)\,\mathcal T(dy\mid x,u,\mu)$.
Then
\begin{align*}
\langle \mu_{t+1}^N,f\rangle
=
\frac{1}{N}\sum_{i=1}^N f(x_{t+1}^i).
\end{align*}
Conditioning on \(\{(x_t^i,u_t^i)\}_{i=1}^N\), we obtain
\begin{align*}
\mathbb E\!\left[\langle \mu_{t+1}^N,f\rangle \,\middle|\, \{(x_t^i,u_t^i)\}_{i=1}^N\right]
&=
\frac{1}{N}\sum_{i=1}^N \mathcal T f(x_t^i,u_t^i,\mu_t^N)\\
&=
\langle \theta_t^N,\mathcal T f(\cdot,\cdot,\mu_t^N)\rangle.
\end{align*}
By our assumptions, for every \(f\in C_b(\mathbb X)\), the map $(x,u,\mu)\mapsto \mathcal T f(x,u,\mu)$ is bounded and continuous on \(\mathbb X\times\mathbb U\times\mathcal P(\mathbb X)\).
Hence,
\begin{align*}
\lim_{N\to \infty}\langle \theta_t^N,\mathcal T f(\cdot,\cdot,\mu_t^N)\rangle
=
\langle \theta_t,\mathcal T f(\cdot,\cdot,\mu_t)\rangle
=
\langle F(\mu_t,\theta_t),f\rangle.
\end{align*}
Hence
\begin{align*}
\lim_{N\to \infty}\mathbb E\!\left[\langle \mu_{t+1}^N,f\rangle \,\middle|\, \{(x_t^i,u_t^i)\}_{i=1}^N\right]
=
\langle F(\mu_t,\theta_t),f\rangle.
\end{align*}
Since \(x_{t+1}^{1:n}\) are conditionally independent given \(\{(x_t^i,u_t^i)\}_{i=1}^N\),
\begin{align*}
\var\!\left(\langle \mu_{t+1}^N,f\rangle \,\middle|\, \{(x_t^i,u_t^i)\}_{i=1}^N\right)
&=
\frac{1}{N^2}\sum_{i=1}^N
\var\!\left(f(x_{t+1}^i)\mid x_t^i,u_t^i,\mu_t^N\right)\\
&\le \frac{4\|f\|_\infty^2}{N}.
\end{align*}
Therefore, by Chebyshev's inequality, for every $\epsilon>0$
\begin{align*}
\mathbb P\!\left(
\left|
\langle \mu_{t+1}^N,f\rangle
-
\mathbb E\!\left[\langle \mu_{t+1}^N,f\rangle \,\middle|\, \{(x_t^i,u_t^i)\}_{i=1}^N\right]
\right|>\epsilon
\,\middle|\,
\{(x_t^i,u_t^i)\}_{i=1}^N
\right)
\le
\frac{4\|f\|_\infty^2}{N\epsilon^2}
\to 0
\end{align*}
as $N\to \infty$. This yields that $\langle \mu_{t+1}^N,f\rangle$ converges to 
$\langle F(\mu_t,\theta_t),f\rangle$
in probability for every \(f\in\{f_m\}\) as $N\to \infty$. Since \(\rho\) metrizes the weak topology on \(\mathcal P(\mathbb X)\), we have
\begin{align*}
&\mathbb P\!\left(
\rho(\mu_{t+1}^N,F(\mu_t,\theta_t))>\epsilon
\,\middle|\,
\{(x_t^i,u_t^i)\}_{i=1}^N
\right)\\
&=
\mathbb P\!\left(
\sum_{m=1}^{\infty} 2^{-(m+1)}
\left|
\langle \mu_{t+1}^N,f_m\rangle-\langle F(\mu_t,\theta_t),f_m\rangle
\right|
>\epsilon
\,\middle|\,
\{(x_t^i,u_t^i)\}_{i=1}^N
\right) \\
&\le
\mathbb P\!\left(
\sum_{m=1}^{k} 2^{-(m+1)}
\left|
\langle \mu_{t+1}^N,f_m\rangle-\langle F(\mu_t,\theta_t),f_m\rangle
\right|
>\epsilon/2
\,\middle|\,
\{(x_t^i,u_t^i)\}_{i=1}^N
\right) \\
&\quad+
\mathbb P\!\left(
\sum_{m=k+1}^{\infty} 2^{-(m+1)}
\left|
\langle \mu_{t+1}^N,f_m\rangle-\langle F(\mu_t,\theta_t),f_m\rangle
\right|
>\epsilon/2
\,\middle|\,
\{(x_t^i,u_t^i)\}_{i=1}^N
\right).
\end{align*}
Since
$\|f_m\|_\infty\le 1$, we can choose \(k\) large enough so that $\sum_{m=k+1}^{\infty} 2^{-(m+1)} 2 < \epsilon/2$.
Then for this \(k\),
\begin{align*}
\lim_{N\to\infty}\mathbb P\!\left(
\sum_{m=1}^{k} 2^{-(m+1)}
\left|
\langle \mu_{t+1}^N,f_m\rangle-\langle F(\mu_t,\theta_t),f_m\rangle
\right|
>\epsilon/2
\,\middle|\,
\{(x_t^i,u_t^i)\}_{i=1}^N
\right)
= 0.
\end{align*}

Hence, conditional on \(\{(x_t^i,u_t^i)\}_{i=1}^N\),
$\mu_{t+1}^N$ converges to $F(\mu_t,\theta_t)$
in probability in $\mathcal{P}(\mathbb X)$ as $N\to \infty$. This concludes that $\mathcal{L}(\mu_{t+1}^n \mid \mu_t^n,\theta_t^n)(\cdot)$ converges to 
$\delta_{F(\mu_t,\theta_t)}(\cdot)$ as $n\to \infty$.

Hence, by \cite{kara2020robustness} or \cite[Theorem 5.1]{langen1981convergence}, the value function $J_{0,T}^{n}$ of the finite population measure-valued MDP converges to that of the limiting MDP $J_{0,T}$.

\item {\bf Step 3: Convergence of empirical measures under symmetric and independent policies.}
Consider the limiting strategic measure $P_{\pi}^{\star}$ induced by symmetric and independent randomized policy $(\pmb{\pi^{\star}},\pmb{\pi^{\star}}, \ldots)$. Under such a policy, the induced states $(x^1_{t}, x_{t}^2, \ldots)$ is infinite exchangeable for every $t\geq 0$. Denote their directing measure by $\mu_t$ (which might be a random measure on $\mathcal{P}(\X)$). Since they are independent conditioned on $\mu_t$, the empirical measure $\mu_t^N$ converges weakly almost surely to $\mu_t$ as $N\to \infty$. Since the policies are symmetric and independent, $((x^1_{t}, u_{t}^1), (x^2_{t}, u_{t}^2), \ldots)$ are also infinitely-exchangeable. Hence, their empirical measure $\theta_t^N$ converges weakly to their directing measure of $\theta_t$ as $N\to \infty$. From the MDP formulation, we have $\mu_{t+1}^N \sim \eta^N(\cdot |\mu_t^N, \theta_t^N)$, which, following from our analysis in Step 2, converges to $\mu_{t+1} \sim \eta(\cdot|\mu_t,\theta_t)$.

Also, the running cost becomes
\begin{align*}
\widetilde{c}(\theta_t,\mu_{t})=\int {c}(x_{t}^{R},u_{t}^{R},\mu_{t}^{R}) \theta_t(dx_{t}^{R},du_{t}^{R}).
\end{align*}
Hence, for symmetric independent policies, we can write the following value iterations for the representative agent  
    \begin{equation*}
\left\{
\begin{aligned}
J_{T-1,T}^R(\mu^R)&:=\inf_{\theta^R\in U(\mu^R)}\widetilde{c}(\mu^R, \theta^R),\\
J_{t,T}^R(\mu^R)
&:=\inf_{\theta^R\in U(\mu^R)}\left\{\widetilde{c}(\mu^R, \theta^R)+\beta \int J_{t+1,T}^R(\mu_{+}) \eta(d\mu_{+}\mid\mu^R, \theta^R)\right\}.
\end{aligned}
\right.
\end{equation*}

\item {\bf Step 4: Optimality of the limiting symmetric and independent policy.} 
Comparing $J_{0,T}$ and $J_{0,T}^{R}$, implies that
$$\inf_{{P_{\pi}}\in L} J_{T}({P_{\pi}})\leq \inf_{{P_{\pi}}\in \LPRS} J_{T}({P_{\pi}})=J_{0,T}^{R}\leq J_{0,T},$$ Hence, we can conclude that $$\inf_{{P_{\pi}}\in L} J_{T}({P_{\pi}})= J_{0,T},$$ and hence,  $\pmb{\theta^{\star}_t}$ with marginals $\pmb{\mu^{\star}_t}$ which are induced by the limiting strategic measure ${P_{\pi}^{\star}}$ are optimal for \PIT. This implies that ${P_{\pi}^{\star}}$ is optimal for \PIT. Since ${P_{\pi}^{\star}}\in \LEX$, we have ${P_{\pi}^{\star}}\in \LCOS$. Hence, $(\underline{\pmb{{x}^{\star}}},\underline{\pmb{{u}^{\star}}})$ are i.i.d. conditioned on common randomness, which is its directing measure $\pmb{\theta^{\star}}$ (with marginal on states as $\pmb{\mu^{\star}}$). {By Theorem \ref{the:ald}, $\pmb{\theta^{\star}}$, the limit in distribution of the empirical measure $\pmb{\theta^{n\star}}$, is a directing 
measure of $(\underline{\pmb{{x}^{\star}}},\underline{\pmb{{u}^{\star}}})$, and hence, the common randomness $z$ corresponds to the directing measure $\pmb{\theta^{\star}}$ and can be computed by the agents as $n\to \infty$.}  Additionally, the following value iterations randomized policies in ${P_{\pi}^{\star}}$ are Markovian.

\item {\bf Step 5: Limit as $T\to \infty$ for infinite horizon problem.} By the monotone convergence theorem, the limit of $J^{\star}_{t,T}$ exists as $T\to \infty$ using a fixed point argument, and hence, we conclude part (i) by establishing a similar result for \PIN\ and showing that randomized policies are additionally stationary.}
\end{itemize}

Part (ii): Following part (i), there exists an optimal strategic measure that is conditionally i.i.d. given the weak limit $\pmb{\theta}$ (with marginal $\pmb{\mu}$) of the empirical measure on states and actions for \PIT\ and \PIN. Hence, the result follows from \eqref{eq:rep-policy} since the limiting policy can be realized by a symmetric and independent policy under MF information sharing. 
\end{proof}

\section{Proof of Section \ref{sec:verif-0}.}

\subsection{Proof of Corollary \ref{Cor:MV}.}\label{APP:CorMV}

Following Theorem \ref{the:3}, we can, without loss, study the MF limit by a symmetric and independent policy under MF information sharing. By the ergodic theorem, we have the empirical measures $\mu^{N}_0$ converge weakly almost surely to $\mu_0=\mathcal{L}(x_0^{R})$ as $N\to \infty$. Since policies are symmetric and independent, the induced state-action pairs $((x_0^1,u_0^1),(x_0^2,u_0^2), \ldots)$ are infinitely exchangeable, and hence, by the ergodic theorem, the empirical measure $\theta_0^N$ converges weakly almost surely to the directing measure $\theta_0$ with marginal $\mu_0$, conditional on $\theta_0$. Since the policy is selected in the limit using $\mu_0$ (which is deterministic), we conclude that the induced probability measure on the action with a deterministic policy is also deterministic. Hence, $\theta_0$ is deterministic and $\theta_0=\mathcal{L}(x_0^R,u_0^R)$. Hence, we conclude that the empirical measure $\theta^N_0$ converges weakly almost surely to $\theta_0=\mathcal{L}(x_0^R,u_0^R)$. Following from an analogous argument as that in step 2 of the proof of Theorem \ref{the:3}, we get $\mu_{1}^{N}$ converges weakly to $\mu_1=\mathcal{L}(x_1^R|\theta_0)$.

Since $\mu_0$ and $\theta_0$ are deterministic, we conclude that $\mu_{1}=\mathcal{L}(x_1^R)$. Utilizing the same argument sequentially, we conclude that the dynamics describing the propagation of $\mu_{t+1}$ given $\mu_t$ and $\theta_t$ is deterministic since the induced probability measure on $\mu_{t+1}$ by kernel $\eta$ is atomic. Hence, the value iterations corresponding to that in Step 3 of Theorem \ref{the:3} and the representative agent are equivalent. This implies optimality of the McKean-Vlasov representative control problem and completes the proof.

\section{Proof of Section \ref{sec:approximation}.}

\subsection{Proof of Theorem \ref{the:4}.}\label{APP:the4}

Following Theorem \ref{the:2}, we can restrict our search for an optimal strategic measure to $\LEXN$ for \PNT\ without any loss of optimality. Let $P_{T}^{N\star} \in  \LEXN$ be an optimal strategic measure for \PNT, i.e.,
\begin{align}\label{eq:asy32.1}
   J^{N}_{T}(P_{T}^{ N\star}) \leq \inf\limits_{P_{T}^{N} \in \LEXN}J^{N}_{T}(P_{T}^{ N})= J_{0,T}^N. 
\end{align}
Following from the proof of Theorem \ref{the:3}, using $\{P_{T}^{N\star}\}_{N} \subseteq  \LEXN$, by considering the indexes as a sequence of i.i.d. random variables with uniform distribution on the set $\{1,\dots,N\}$, we can construct a sequence of  infinitely exchangeable policies $\{P_{T}^{N,\infty\star}\}$ where the restriction of an infinitely exchangeable policy to $N$ first components, $P^{\infty\star}_{T, N} \in \LEX\big{|}_{N}$, satisfies 
\begin{flalign}
& J^{N}_{T}(P^{\infty\star}_{T, N}) \leq  J^{N}_{T}(P_{T}^{ N\star})+\widehat{\epsilon_{N}
}\label{eq:asymendif}
\end{flalign}
for some $\widehat{\epsilon_{N}
}\geq 0$ such that $\widehat{\epsilon_{N}
}$ goes to $0$ as $N\to \infty$. Hence, \eqref{eq:asy32.1} and \eqref{eq:asymendif} imply that
\begin{flalign*}
\inf\limits_{P_{T}^{N} \in \LEX|_N}J^{T}_{N}(P_{T}^{ N})\leq \inf\limits_{P_{T}^{N} \in \LEXN}J^{T}_{N}(P_{T}^{ N})+\widehat{\epsilon_{N}
}.
\end{flalign*}
 By de Finetti's theorem, $\LEX|_N=\LCONS$. Following from the proof of Theorem \ref{the:3}, we additionally get
 \begin{flalign}\label{eq:approx-exc}
\inf\limits_{P_{T}^{N} \in \LPRS^N}J^{T}_{N}(P_{T}^{ N})\leq \inf\limits_{P_{T}^{N} \in \LEXN}J^{N}_{T}(P_{\pi}^{ N})+\widehat{\epsilon_{N}
}.
\end{flalign}
 Let $P_{T}^{\star}\in \LPRS$ be an optimal policy for \PIT\ and $P_{T,N}^{\star}$ is the restriction of $P_{T}^{
\star}$ to the first $N$ components. Since $P_{T,N}^{\star}$ is symmetric and independent,   following from proof of Theorem \ref{the:3}, we have
\begin{flalign}
J^{N}_{T}(P_{T,N}^{\star})\leq \inf\limits_{P_{T}^{N} \in \LPRNS} J_{T}^{N}(P_{T}^{N})+\widetilde{\epsilon_{N}}\label{eq:L23}
\end{flalign}
for some $\widetilde{\epsilon_{N}}\geq 0$ such that $\widetilde{\epsilon_{N}
}$ goes to $0$ as $N\to \infty$. Letting $\epsilon_{N}=\widetilde{\epsilon_{N}}+\widehat{\epsilon_{N}}$,  \eqref{eq:L23} and \eqref{eq:approx-exc} imply that
\begin{flalign}
J_{T}^{N}(P_{T,N}^{\star})\leq \inf\limits_{P_{\pi}^{N} \in \LEX^N} J_{T}^{N}(P_{T}^{N})+{\epsilon_{N}},\label{eq:approxl}
\end{flalign}
which completes the proof thanks to \eqref{eq:asy32.1}. Similarly, we can show the result for \PN\ using the optimal strategic measure of \PIN.

\subsection{Proof of Theorem \ref{the:Approx-MK}.}\label{APP:the5}
Since $\mu_{0}$ be non-atomic and $\mathcal{T}(\cdot|x_0^R,u_{0}^{R},\mu_{0})$ is non-atomic for every $x_0^R,u_{0}^{R}$ and $\mu_{0}$, we have $\{\mu_t\}_{t\geq0}$ is non-atomic under any policy. We first prove the result for the finite horizon problems. We start from $T-1$. Following \cite{saldi2018finite,saldi2017asymptotic} or  \cite{kara2023q}, for every $\widehat{\mu}\in \mathcal{P}(\widehat{\mathbb{X}})$, we have
\begin{align*}
   |\widehat{J}_{T-1,T}^n(\widehat{\mu})-{J}_{T-1,T}(\widehat{\mu})|
   &=\left|\inf_{\widehat\theta\in \widehat{U}_n(\widehat{\mu})}\widehat{c}(\widehat{\mu}, \widehat\theta)-\inf_{\theta\in {U}(\widehat{\mu})}\widehat{c}(\widehat{\mu}, \theta)\right|\\
   &=\bigg|\int c(x,u,\widehat{\mu}) \pi^{\star}_n(du|x)\widehat{\mu}(dx) -\int c(x,u,\widehat{\mu}) \pi^{\star}(du|x)\widehat{\mu}(dx)\bigg|\to 0
\end{align*}
where $\pi^{\star}_n$ is an optimal policy at $T-1$ for the quantized ${\sf MDP}_{n}$ and $\pi^{\star}$ is an optimal policy at $T-1$ for the original problem, fixing $\widehat{\mu}$.
Suppose that $\widehat{\mu}_{n}$ converges weakly to $\mu$ as $n\to \infty$. We also have
\begin{align*}
   &|\widehat{J}_{T-1,T}^n(\widehat{\mu}_n)-{J}_{T-1,T}({\mu})|\\
   &=\left|\inf_{\widehat\theta\in \widehat{U}_n(\widehat{\mu}_n)}\widehat{c}(\widehat{\mu}_n, \widehat\theta)-\inf_{\theta\in {U}({\mu})}\widehat{c}({\mu}, \theta)\right|\\
   &=\bigg|\int c(x,u,\widehat{\mu}_n) \pi^{\star}_n(du|x)\widehat{\mu}_n(dx) -\int c(x,u,{\mu}) \pi^{\star}(du|x){\mu}(dx)\bigg|\to 0.
\end{align*}
This implies that $\widehat{J}_{T-1,T}^n$ converges continuously to ${J}_{T-1,T}$. Next, we let $t=T-2$. We have 
\begin{align*}
    &|\widehat{J}_{T-2,T}^n(\widehat{\mu})-{J}_{T-2,T}(\widehat{\mu})|\\
    &=\bigg{|}\inf_{\widehat\theta\in \widehat{U}_n(\widehat{\mu})}\left\{\widehat{c}(\widehat{\mu}, \widehat\theta)+\beta\int \widehat{J}_{T-1,T}^n(\widehat{\mu}_{T-1}) \widehat{\eta}_n(d\widehat{\mu}_{T-1}|\widehat{\mu}_{T-2},\widehat\theta)\right\}\\
    &-\inf_{\theta\in {U}(\widehat{\mu})}\left\{\widehat{c}(\widehat{\mu}, \theta)+\beta\int {J}_{T-1,T}({\mu}_{T-1}) {\eta}(d{\mu}_{T-1}|\widehat{\mu}_{T-2},\widehat\theta)\right\}\bigg{|}\\
    &=\left|\int c(x,u,\widehat{\mu}) \pi^{\star}_{n}(du|x)\widehat{\mu}(dx)-\int c(x,u,\widehat{\mu}) \pi^{\star}(du|x)\widehat{\mu}(dx)\right|\\
    &+\beta\left|\int \widehat{J}_{T-1,T}^n(\widehat{\mu}_{T-1}) \widehat{\eta}_n(d\widehat{\mu}_{T-1}|\widehat{\mu}_{T-2},\widehat\theta^{\star}_n)-\int {J}_{T-1,T}({\mu}_{T-1}) {\eta}(d{\mu}_{T-1}|\widehat{\mu}_{T-2},\theta^{\star})\right|
\end{align*}
The first term above goes to zero as $n\to \infty$. We have
\begin{align*}
    &\left|\int \widehat{J}_{T-1,T}^n(\widehat{\mu}_{T-1}) \widehat{\eta}_n(d\widehat{\mu}_{T-1}|\widehat{\mu}_{T-2},\widehat\theta^{\star}_n)-\int {J}_{T-1,T}({\mu}_{T-1}) {\eta}(d{\mu}_{T-1}|\widehat{\mu}_{T-2},\theta^{\star})\right|\\
    &\leq \left|\int \widehat{J}_{T-1,T}^n(\widehat{\mu}_{T-1}) \widehat{\eta}_n(d\widehat{\mu}_{T-1}|\widehat{\mu}_{T-2},\widehat\theta^{\star}_n)-\int {J}_{T-1,T}({\mu}_{T-1}) {\eta}(d{\mu}_{T-1}|\widehat{\mu}_{T-2},\widehat\theta^{\star}_n)\right|\\
    &+\left|\int {J}_{T-1,T}({\mu}_{T-1}) {\eta}(d{\mu}_{T-1}|\widehat{\mu}_{T-2},\widehat\theta^{\star}_n)-\int {J}_{T-1,T}({\mu}_{T-1}) {\eta}(d{\mu}_{T-1}|\widehat{\mu}_{T-2},\theta^{\star})\right|
\end{align*}
Since $\eta$ is weakly continuous in $\theta$, the last term converges to zero. Since $\widehat{J}_{T-1,T}^n$ converges continuously to ${J}_{T-1,T}$, we have
\begin{align*}
    &\left|\int \widehat{J}_{T-1,T}^n(\widehat{\mu}_{T-1})\widehat{\eta}_n(d\widehat{\mu}_{T-1}|\widehat{\mu}_{T-2},\widehat\theta^{\star}_n)-\int {J}_{T-1,T}({\mu}_{T-1}) {\eta}(d{\mu}_{T-1}|\widehat{\mu}_{T-2},\widehat\theta^{\star}_n)\right|\\
    &\leq\sup_{\theta \in {U}(\widehat{\mu})}  \left|\int \widehat{J}_{T-1,T}^n({\mu}_{T-1})\widehat{\eta}_n(d{\mu}_{T-1}|\widehat{\mu}_{T-2},\theta)-\int {J}_{T-1,T}({\mu}_{T-1}) {\eta}(d{\mu}_{T-1}|\widehat{\mu}_{T-2},\theta)\right|\to 0.
\end{align*}
Following an analogous argument used above, we establish \eqref{the:Ap-main}.

We now establish the result for the infinite horizon. Let $\epsilon>0$. We have
\begin{align*}
    |J_{\infty}^{n}(\mu_0)-J_{\infty}(\mu_0)|
    &\leq |J_{\infty}^{n}(\mu_0)-J_{0,T}^{n}(\mu_0)|+|J_{0,T}^{n}(\mu_0)-J_{0,T}(\mu_0)|+|J_{0,T}^{n}(\mu_0)-J_{0,T}^{n}(\mu_0)|\\
    &\leq \epsilon_{T} + \epsilon_{T}^n. 
\end{align*}
Since $\epsilon_{T}$ converges to zero as $T\to \infty$, and $\epsilon_{T}^n$ converges to zero as $n$ goes to infinity, for $n$ large enough, we can select $T$ large enough such that $\epsilon_{T} + \epsilon_{T}^n\leq \epsilon$. Since $\epsilon> 0$ is arbitrary, this completes the proof.

\section*{Acknowledgements.}
The research of the first and third authors is supported by the Natural Sciences and Engineering Research Council (NSERC) of Canada.

\bibliographystyle{plain}

\end{document}